\documentclass[12pt]{article}
\usepackage[utf8]{inputenc}
\usepackage[left=1in,right=1in, top=1.2in,bottom=1.2in,includehead]{geometry}
\usepackage{amssymb,amsthm,amsmath,pgf}
\usepackage{enumitem}
\usepackage{mathrsfs}
\usepackage{placeins}
\usepackage{wrapfig}
\usepackage{indentfirst}
\usepackage{fancyhdr}
\usepackage{accents}
\usepackage{algorithm,algpseudocode}
\usepackage{url}
\usepackage[numbers]{natbib}
\usepackage[hidelinks]{hyperref}

\usepackage{lineno}
\usepackage{xcolor}

\newtheorem{theorem}{Theorem}[section]

\newtheorem{corollary}{Corollary}[section]
\newtheorem{conjecture}{Conjecture}[section]

\newtheorem{definition}{Definition}[section]
\newtheorem{notation}{Notation}[section]
\newtheorem{claim}{Claim}

\newcommand{\newcase}[1]{\vspace{0.7em}\noindent\textbf{Case #1: }}
\newcommand{\cbeginproof}[0]{\par\noindent\textit{Proof.} }
\newcommand{\cendproof}[0]{ \qed\par\vspace{1em}}

\title{Fault-tolerant Locating-Dominating Sets \\ on the Infinite King Grid}
\author{
    \small Devin C. Jean\\
    \small Computer Science Department \\
    \small Vanderbilt University\\
    \small \texttt{devin.c.jean@vanderbilt.edu}
    \and
    \small Suk J. Seo\\
    \small Computer Science Department\\
    \small Middle Tennessee State University\\
    \small \texttt{Suk.Seo@mtsu.edu}
}
\date{}

\begin{document}
\maketitle
\thispagestyle{empty}

\begin{abstract}
Let $G$ be a graph of a network system with vertices, $V(G)$, representing physical locations and edges, $E(G)$, representing informational connectivity.
A \emph{locating-dominating (LD)} set $S \subseteq V(G)$ is a subset of vertices representing detectors capable of sensing an ``intruder'' at precisely their location or somewhere in their open-neighborhood---an LD set must be capable of locating an intruder anywhere in the graph.
We explore three types of fault-tolerant LD sets: \emph{redundant LD} sets, which allow a detector to be removed, \emph{error-detecting LD} sets, which allow at most one false negative, and \emph{error-correcting LD} sets, which allow at most one error (false positive or negative).
In particular, we determine lower and upper bounds for the minimum density of these three fault-tolerant locating-dominating sets in the \emph{infinite king grid}.
\end{abstract}

\vspace{1em}


\section{Introduction}

Let $G$ be a graph---with vertex set $V(G)$ and edge set $E(G)$---which models some system or facility with vertices representing physical locations or pieces of hardware and edges representing their physical or informational connectivity.
There may be some ``intruder'' or other undesirable event in the system, which must be automatically detected and located anywhere in the graph via detectors placed at some subset $S \subseteq V(G)$ of vertices.

A detector vertex $v \in S$ is associated with a set of \emph{detection regions} $R(v) \subseteq \mathscr{P}(V(G))$, representing the areas in which it can sense an intruder.
A \emph{detection system} is a set of detector vertices, $S \subseteq V(G)$, such that an intruder can be precisely located anywhere in the graph provided the detection signal (presence or absence of an intruder) for every detection region.
In this paper, we will explore \emph{locating-dominating (LD) sets}, which are a class of detection system whose detectors have two detection regions: $R(v) = \{\{v\}, N(v)\}$, where $N(v)$ denotes the open-neighborhood of $v$: $\{u \in V(G) : uv \in E(G)\}$.
Other choices for $R(v)$ have been explored, such as $R(v) = \{N(v)\}$, yielding \emph{open-locating-dominating (OLD) sets} \cite{kinc15, seo10}, or $R(v) = \{N[v]\}$---where $N[v]$ denotes the closed neighborhood of $v$: $N(v) \cup \{v\}$---yielding \emph{identifying codes (ICs)} \cite{karp98a}.
Jean and Lobstein \cite{dombib} maintain a bibliography containing over 500 papers on these topics and related graphical parameters.

Conceptually, every LD detector vertex $v \in S$ will transmit a value of 0 (denoting no intruder), 1 (denoting an intruder in the open-neighborhood), or 2 (denoting an intruder at $v$).
The goal is then to select a sufficient subset of vertices $S \subseteq V(G)$ so that we can locate an intruder anywhere in the graph.
We are interested in going one step further: we would like the system to be able to function even in the event of an error, known as a \emph{fault-tolerant LD (FTLD) set}.
In particular, we consider three types of FTLD sets: \emph{redundant LD (RED:LD) sets}, in which the system must be fully functioning even if at most one detector is removed or disabled, \emph{error-detecting LD (DET:LD)}, where the system allows at most one false-negative (i.e. sending 0 instead of 1 or 2), and \emph{error-correcting LD (ERR:LD) sets}, where the system allows any single transmission error (false-positive or false-negative).

The following theorems characterize the conditions of LD, RED:LD, DET:LD, and ERR:LD sets.

\begin{theorem}[\cite{jean23b}]\label{theo:ld-char}
A set $S \subseteq V(G)$ is an LD set if and only if the following are true:
\begin{enumerate}[noitemsep, label=\roman*.]
    \item $\forall v \in V(G)-S$, $|N(v) \cap S| \ge 1$
    \item $\forall v,u \in V(G)-S$ with $v \neq u$, $|(N(v) \cap S) \triangle (N(u) \cap S)| \ge 1$
\end{enumerate}
\end{theorem}

\begin{theorem}[\cite{jean23b}]\label{theo:red-ld-char}
A set $S \subseteq V(G)$ is a RED:LD set if and only if the following are true:
\begin{enumerate}[noitemsep, label=\roman*.]
    \item $\forall v \in V(G)$, $|N[v] \cap S| \ge 2$
    \item $\forall v \in S$ and $\forall u \in V(G)-S$, $|((N(v) \cap S) \triangle (N(u) \cap S)) - \{v\}| \ge 1$
    \item $\forall v,u \in V(G)-S$ with $u \neq v$, $|(N(v) \cap S) \triangle (N(u) \cap S)| \ge 2$
\end{enumerate}
\end{theorem}

\begin{theorem}[\cite{jean23a}]\label{theo:det-ld-char}
A set $S \subseteq V(G)$ is a DET:LD set if and only if the following are true:
\begin{enumerate}[noitemsep, label=\roman*.]
    \item $\forall v \in V(G)$, $|N[v] \cap S| \ge 2$
    \item $\forall v,u \in S$ with $u \neq v$, $|(N(v) \cap S) \triangle (N(u) \cap S)| \ge 1$.
    \item $\forall v \in V(G)-S$ and $\forall u \in S$, $|(N(v) \cap S) - (N(u) \cap S)| \ge 2$ or $|(N(u) \cap S) - (N(v) \cap S)| \ge 1$
    \item $\forall v,u \in V(G)-S$ with $u \!\neq\! v$, $\!|(N(v) \cap S) - (N(u) \cap S)| \ge 2$ or $|(N(u) \cap S) - (N(v) \cap S)| \ge 2$
\end{enumerate}
\end{theorem}

\newpage
\begin{theorem}[\cite{errld}]\label{theo:err-ld-char}
A set $S \subseteq V(G)$ is an ERR:LD set if and only if the following are true:
\begin{enumerate}[noitemsep, label=\roman*.]
    \item $\forall v \in V(G)$, $|N[v] \cap S| \ge 3$
    \item $\forall v,u \in S$ with $u \neq v$, $|((N(v) \cap S) \triangle (N(u) \cap S)) - \{v, u\}| \ge 1$
    \item $\forall v \in V(G)-S$ and $\forall u \in S$, $|((N(v) \cap S) \triangle (N(u) \cap S)) - \{u\}| \ge 2$
    \item $\forall v,u \in V(G)-S$ with $u \neq v$, $|(N(v) \cap S) \triangle (N(u) \cap S)| \ge 3$
\end{enumerate}
\end{theorem}

For an LD-based detection system, a vertex $x$ is said to be \emph{$k$-dominated} if $|N[v] \cap S| \ge k$, and a vertex pair $(v,u)$ is said to be \emph{$k$-distinguished} if $|((N(v) \cap S) \triangle (N(u) \cap S)) - \{v,u\}| \ge k$.
When discussing any particular type of detection system with a known characterization, we will also say that a vertex pair is simply ``distinguished'' if it satisfies said characterization.
The properties given in Theorems \ref{theo:red-ld-char} and \ref{theo:err-ld-char} can be rewritten as the following, which make use of $k$-domination and $k$-distinguishing.

\begin{corollary}[\cite{jean23b}]\label{cor:red-ld}
A set $S \subseteq V(G)$ is a RED:LD set if and only if the following are true:
\begin{enumerate}[noitemsep, label=\roman*.]
    \item $\forall v \in V(G)$, $|N[v] \cap S| \ge 2$
    \item $\forall v,u \in V(G)$ with $u \neq v$, $|((N(v) \cap S) \triangle (N(u) \cap S)) - \{v, u\}| \ge 2 - |\{v,u\} \cap S|$
\end{enumerate}
\end{corollary}

\begin{corollary}[\cite{errld}]\label{cor:err-ld}
A set $S \subseteq V(G)$ is an ERR:LD set if and only if the following are true:
\begin{enumerate}[noitemsep, label=\roman*.]
    \item $\forall v \in V(G)$, $|N[v] \cap S| \ge 3$
    \item $\forall v,u \in V(G)$ with $u \neq v$, $|((N(v) \cap S) \triangle (N(u) \cap S)) - \{v, u\}| \ge 3 - |\{v,u\} \cap S|$
\end{enumerate}
\end{corollary}

Given a graph, $G$, we are interested in the smallest possible sets with the required properties.
Jean and Seo \cite{jean23b, jean23a, errld} have studied the optimal sets for special classes of graphs such as trees, cubic graphs, and infinite grids.
They also proved that these minimization problems on arbitrary graphs are NP Complete.
In this paper, we consider optimal sets for these three fault-tolerant LD set parameters in the \emph{infinite king grid (K)}, an 8-regular graph inspired by a chess board, with edges denoting the typical eight moves of a king.
In particular, we explore upper and lower bounds for the smallest RED:LD, DET:LD, and ERR:LD sets on $\textrm{K}$.

The \emph{density} of a subset $S \subseteq V(G)$ is the ratio of $|S|$ to $|V(G)|$, or more formally, $\limsup_{r \to \infty} \frac{|B_r(v) \cap S|}{|B_r(v)|}$ for some center point $v \in V(G)$ where $B_r(v) = \{u \in V(G) : d(u,v) \le r\}$ is the ball of radius $r$ about $v$.
In this paper, we will only explore graphs which satisfy the ``slow-growth'' property \cite{samp24} of $\lim_{r \to \infty} \frac{|B_{r+1}(v)|}{|B_r(v)|} = 1$, which guarantees that the density is invariant of the center point $v$, and that the density of a periodic tiling is equal to the density within any single tiling.
The notations $\textrm{RED:LD}(G)$, $\textrm{DET:LD}(G)$, and $\textrm{ERR:LD}(G)$ denote the minimum density of such a set on $G$.

In Section~\ref{sec:prev-results}, we discuss some of the previous results of bounds for optimal (minimum) fault-tolerant variants of LD, OLD, and IC sets.
Section~\ref{sec:lower-bound-techniques} introduces the concepts of ``share arguments'' and ``discharging techniques,'' which are used later in Section~\ref{sec:lower-bounds} to determine bounds for RED:LD, DET:LD, and ERR:LD sets on the infinite king grid.

\section{Previous Results}\label{sec:prev-results}

The concept of LD sets was originally introduced by Slater \cite{slat87a, slat88a}.
Slater \cite{slat02a} also explored ``FTLD'' sets (equivalent to DET:LD sets in this paper), and proved upper and lower bounds for optimal DET:LD sets in the \emph{infinite square grid (SQ)}: $\left[\frac{12}{23}, \frac{3}{5}\right]$.
Jean and Seo \cite{jean23b, errld} introduced two new fault-tolerant LD sets, RED:LD and ERR:LD sets, and fully characterized and established existence criteria for RED:LD, DET:LD, and ERR:LD sets; these characterizations are included in this paper as Theorems \ref{theo:red-ld-char}, \ref{theo:det-ld-char}, and \ref{theo:err-ld-char}.
Note that characterizations of general fault-tolerant variants of distinguishing sets were established by Seo and Slater \cite{seo15}, though this does not cover LD set variants due to having more than one detection region.
Jean and Seo \cite{jean23a} also proved upper and lower bounds for optimal RED:LD, DET:LD, and ERR:LD sets in SQ: $\left[\frac{2}{5}, \frac{7}{16}\right]$ for redundant, $\left[\frac{12}{23}, \frac{7}{12}\right]$ for error-detecting (improved from \cite{slat02a}), and $\left[\frac{24}{37}, \frac{2}{3}\right]$ for error-correcting.

A variety of graphical parameters have been studied on the infinite king grid \cite{dant18b, fouc13, honk06c, pelt12b}; for instance, Charon et al. \cite{char02a} proved $\textrm{IC}(\textrm{K}) = \frac{2}{9}$, Honkala and Laihonen \cite{honk05c} determined that $\textrm{LD}(\textrm{K}) = \frac{1}{5}$, Seo \cite{seo18} found that $\frac{6}{25} \leq \textrm{OLD}(\textrm{K}) \leq \frac{1}{4}$, and Jean and Seo \cite{ourking} proved $\frac{3}{10} \leq \textrm{RED:OLD}(\textrm{K}) \leq \frac{1}{3}$ and $\frac{3}{11} \leq \textrm{RED:IC}(\textrm{K}) \leq \frac{1}{3}$.

Similar types of fault-tolerant open-locating-dominating sets on other infinite graphs have been studied, as well; for instance on the \emph{infinite hexagonal grid (HEX)} and the \emph{infinite triangular grid (TRI)}.
Seo and Slater \cite{seo15} proved that $\textrm{RED:OLD}(\textrm{HEX}) = \frac{2}{3}$, $\textrm{DET:OLD}(\textrm{HEX}) = \frac{6}{7}$, $\textrm{RED:OLD}(\textrm{SQ}) = \frac{1}{2}$, $\textrm{DET:OLD}(\textrm{SQ}) = \frac{3}{4}$, $\textrm{RED:OLD}(\textrm{TRI}) = \frac{3}{8}$, and $\textrm{DET:OLD}(\textrm{TRI}) \le \frac{5}{9}$.
Later, Jean and Seo \cite{jean23} proved that $\textrm{DET:OLD}(\textrm{TRI})  = \frac{1}{2}$.

\section{Lower Bound Techniques}\label{sec:lower-bound-techniques}

In this section, we explain the concepts of ``share'' and ``discharging,'' which will be used in Section~\ref{sec:lower-bounds} to prove the lower bounds for RED:LD, DET:LD, and ERR:LD sets in the infinite king grid.
We also introduce several notations that will be used to shorten these proofs.

\subsection{Share}\label{sec:share}

Instead of finding a lowerbound for the optimal density directly, it is often easier to find an upperbound for the average ``contribution'' of each detector vertex.
This ``contribution'' is measured by the \emph{share} of a detector vertex, which was first introduced by Slater \cite{slat02a}.
Let $S \subseteq V(G)$ be a set of LD detectors.

\begin{definition}
The \emph{domination count} of a vertex $v \in V(G)$, denoted $dom(v)$, is $|N[v] \cap S|$.
\end{definition}

\begin{definition}
The \emph{partial share} of a vertex $v \in V(G)$, denoted $sh[v]$, is $\frac{1}{dom(v)}$.
\end{definition}

\begin{definition}
The (total) \emph{share} of a detector vertex $v \in S$, denoted $sh(v)$, is $\sum_{u \in D}{sh[u]}$ where $D$ is the set of vertices dominated by $v$.
\end{definition}

For LD detectors, $dom(v) = |N[v] \cap S|$ and $sh(v) = \sum_{u \in N[v]}{sh[u]}$.
If $S$ is a \emph{dominating set}---that is, $\forall v \in V(G)$, $dom(v) \ge 1$---it is clear that, for a finite graph, $\sum_{x \in S}{sh(x)} = |V(G)|$ because each vertex which is $k$-dominated has partial share $\frac{1}{k}$ and this fraction is repeated $k$ times due to having $k$ dominators.
By extension, the average share will be equal to $\frac{|V(G)|}{|S|}$, the inverse density of $S$ in $V(G)$, even in an infinite graph.

\begin{definition}
For some $A \subseteq V(G)$, let $sh[A]$ denote the \emph{block share} of $A$, $\sum_{x \in A}{sh[x]}$.
\end{definition}

For example, $sh[\{v_1,v_2\}] = sh[v_1] + sh[v_2]$, though we often shorten this to $sh[v_1v_2]$.

\begin{notation}
Let $D_k$ denote the set of vertices which are exactly $k$-dominated, $\{v \in V(G) : dom(v) = k\}$, and $D_{k+}$ denote the set of vertices which are at least $k$-dominated, $\cup_{j \ge k}{D_j}$.
\end{notation}

\begin{notation}
Let $\sigma_A$ be $\sum_{x \in A}{\frac{1}{x}}$ where $A$ is a finite sequence of single-character symbols denoting natural numbers (e.g. $\sigma_{123} = \frac{1}{1} + \frac{1}{2} + \frac{1}{3}$ or $\sigma_{abc} = \frac{1}{a} + \frac{1}{b} + \frac{1}{c}$ where $a,b,c \in \mathbb{N}$).
\end{notation}

\begin{wrapfigure}{r}{0.2\textwidth}
    \centering
    \includegraphics[width=0.2\textwidth]{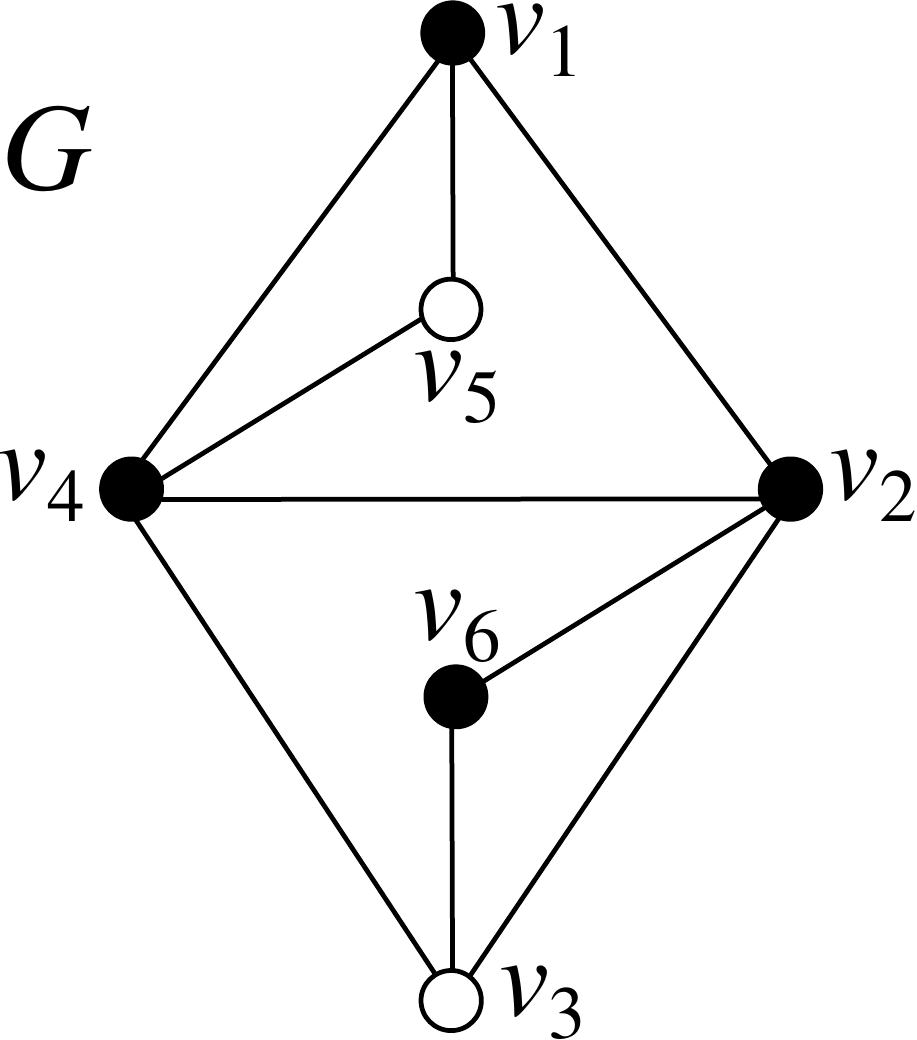}
    \caption{Example RED:LD set.}
    \label{fig:share-demo}
\end{wrapfigure}

Consider the graph $G$ in Figure~\ref{fig:share-demo}, where shaded vertices represent detectors, elements of RED:LD set $S= \{v_1, v_2,v_4,v_6\}$.
We observe that $\{v_5,v_6\} \subseteq D_2$, $\{v_1,v_3,v_4\} \subseteq D_3$, and $v_2 \in D_4$, thus $sh[v_5] = sh[v_6] = \frac{1}{2}$, $sh[v_1] = sh[v_3] = sh[v_4] = \frac{1}{3}$, and $sh[v_2] = \frac{1}{4}$.
Therefore, $sh(v_1) = sh[N[v_1]] = sh[v_1v_2v_4v_5] = \sigma_{3432} = \frac{17}{12}$ and $sh(v_2) = sh[N[v_2]] = sh[v_1v_2v_3v_4v_6] = \sigma_{34332} = \frac{7}{4}$.
Similarly, we can show that $sh(v_4) = \frac{7}{4}$ and $sh(v_6) = \frac{13}{12}$.
And we see that the sum of shares for the dominating set $S$ is $\frac{17}{12} + \frac{7}{4} + \frac{7}{4} + \frac{13}{12} = 6 = |V(G)|$, as expected, and the average is indeed $\frac{6}{4} = \frac{3}{2}$, the inverse density.
Although it is not shown, using $S = V(G)$ in the example graph from Figure~\ref{fig:share-demo} would be an ERR:LD set, and using $S = V(G) - \{v_3\}$ would be a DET:LD set.

\subsection{Discharging}\label{sec:discharge}

In proving a maximum (target) value for the average share of all detector vertices, we would ideally like to prove a maximum value for the share of any individual detector vertex, as this would necessarily be a maximum for the average as well.
However, in doing so we might run into ``problem cases'' where the share of a detector vertex exceeds our target.
In these cases we must prove that said problem case also induces one or more other detector vertices with sufficiently low share that the average will still be bounded by the target.
To prove this bound for the average share, we apply a \emph{discharging argument}, in which we attempt to distribute the excess share of each problematic detector vertex to other detectors with share values below the target.
So long as the final result after discharging ensures that each detector has a share value no more than the target, the bound on the average share is proven.
This operation is sound because the sum of shares, and therefore the average share, remains the same.

\begin{notation}
For some $x \in S$, let $\widehat{sh}(x)$ denote the share of $x$ after discharging.
\end{notation}

\begin{notation}
For some $x \in S$, let $sh_{max}(x)$ denote the maximum value of $sh(x)$ in any sub-configuration.
\end{notation}

When discussing problem cases, we consider having a ``center'' vertex, say $x$, which is assumed to have $sh(x) > t$ where $t$ is the target.
In order for the discharging argument to be valid, we require that $\forall u \in S$, $\widehat{sh}(u) \le t$.
We allow discharging only to neighbors: for any $v \in N(x) \cap S$ with $sh(v) < t$, $v$ can accept at most $t - sh_{max}(v)$ discharge in total.
In particular, we allow $x$ to discharge at most $\frac{1}{k}[t - sh_{max}(v)]$ to $v$, where $k \in \mathbb{N}$ is at least as large as the maximum number of neighbors of $v$ which can simultaneously have shares exceeding $t$, including $x$ itself; as a conservative value, $k \ge |u \in N(v) \cap S : sh_{max}(u) > t|$ can be used.
By this construction, it is clear that this discharging rule is sound, as $\widehat{sh}(v)$ will never exceed $t$.

This strategy differs from other discharging approaches which have been used in that it works backwards: instead of determining the excess share of a problem vertex and dividing it equally to its neighbors, we determine for each neighbor the amount of discharge it can receive from the problem vertex.
Instead of uniform outgoing discharge from the problem case, we have uniform incoming discharge to the non-problem neighbors.
This technique is perhaps more powerful for problem cases with asymmetric shares; that is, where the maximum shares of discharge targets are not equal.

\begin{algorithm}[ht]
\caption{Problem case resolution algorithm}
\label{alg:discharge}
\begin{algorithmic}[1]
\Function{DischargeResolves}{$x, t$} \Comment{we assume $x \in S$ with $sh(x) > t$}
\State{$p \gets sh(x)$}
\ForAll{$v \in N(x) \cap S$} \Comment{for any discharge destination $v$}
    \State{$s \gets sh_{max}(v)$} \Comment{find the highest share $v$ can have}
    \If{$s < t$} \Comment{check if discharge to $v$ is possible}
        \State{$k \gets |\{u \in N(v) \cap S : sh_{max}(u) > t\}|$} \Comment{number of discharge sources for $v$}
        \State{$p \gets p - (t - s) / k$} \Comment{safely discharge from $x$ to $v$}
    \EndIf
\EndFor
\State{\Return{$p \le t$}} \Comment{$p$ is now $\widehat{sh}(x)$; check if we resolved}
\EndFunction
\end{algorithmic}
\end{algorithm}

Algorithm~\ref{alg:discharge} formalizes the logic behind the discharging technique we have described.
Suppose we have some vertex $x \in S$ with $sh(x) > t$; this problem case can be resolved if Algorithm~\ref{alg:discharge} returns $true$, given arguments $x$ and $t$.
Note that it is not strictly necessary to test all $v \in N(x) \cap S$, because as soon as $p \le t$ we are done.

\section{Bounds for Fault-Tolerant LD Sets on K}\label{sec:lower-bounds}

\begin{figure}[ht]
    \centering
    \begin{tabular}{c@{\hskip 4em}c}
        \includegraphics[width=0.375\textwidth]{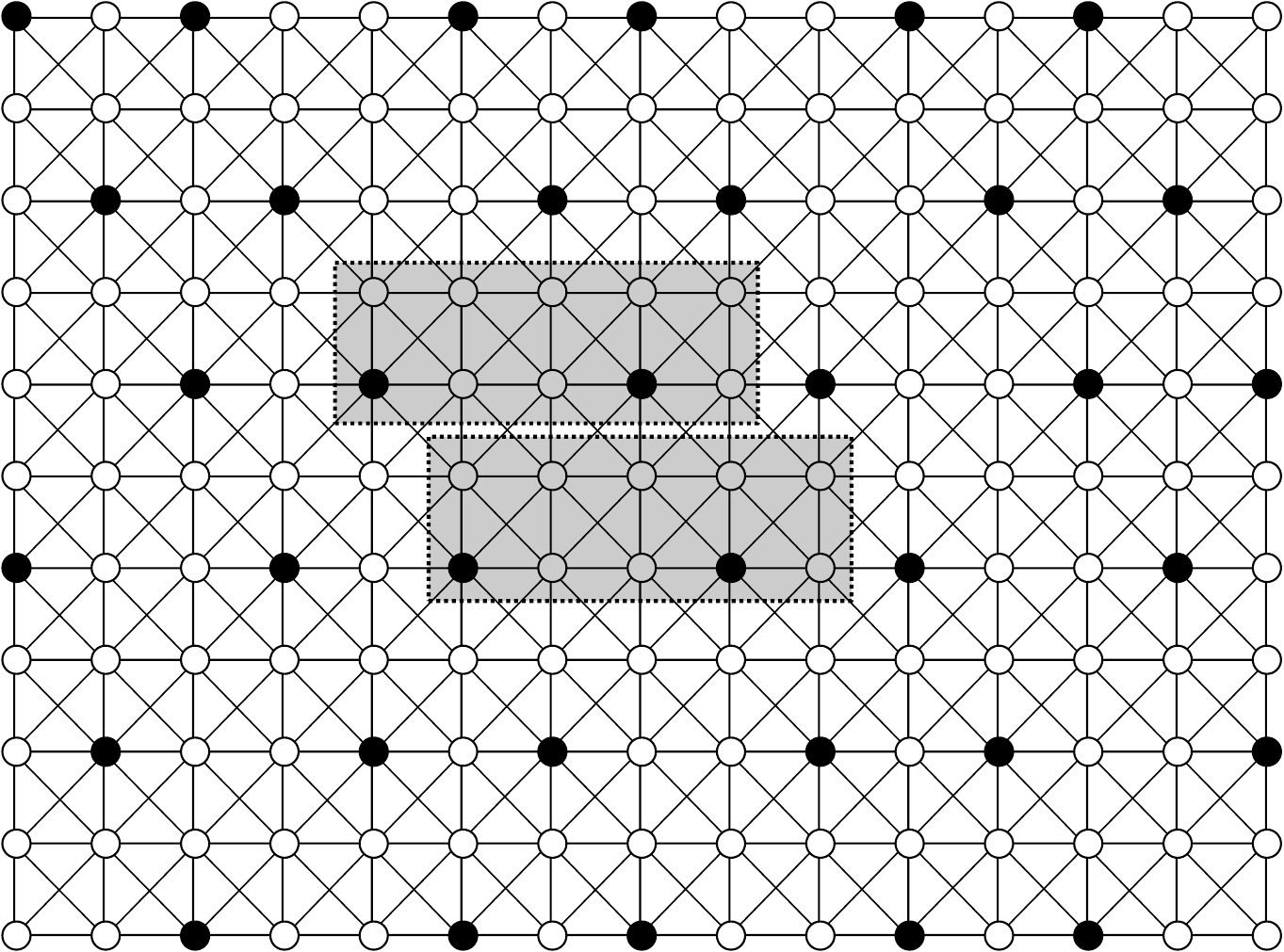} &
        \includegraphics[width=0.375\textwidth]{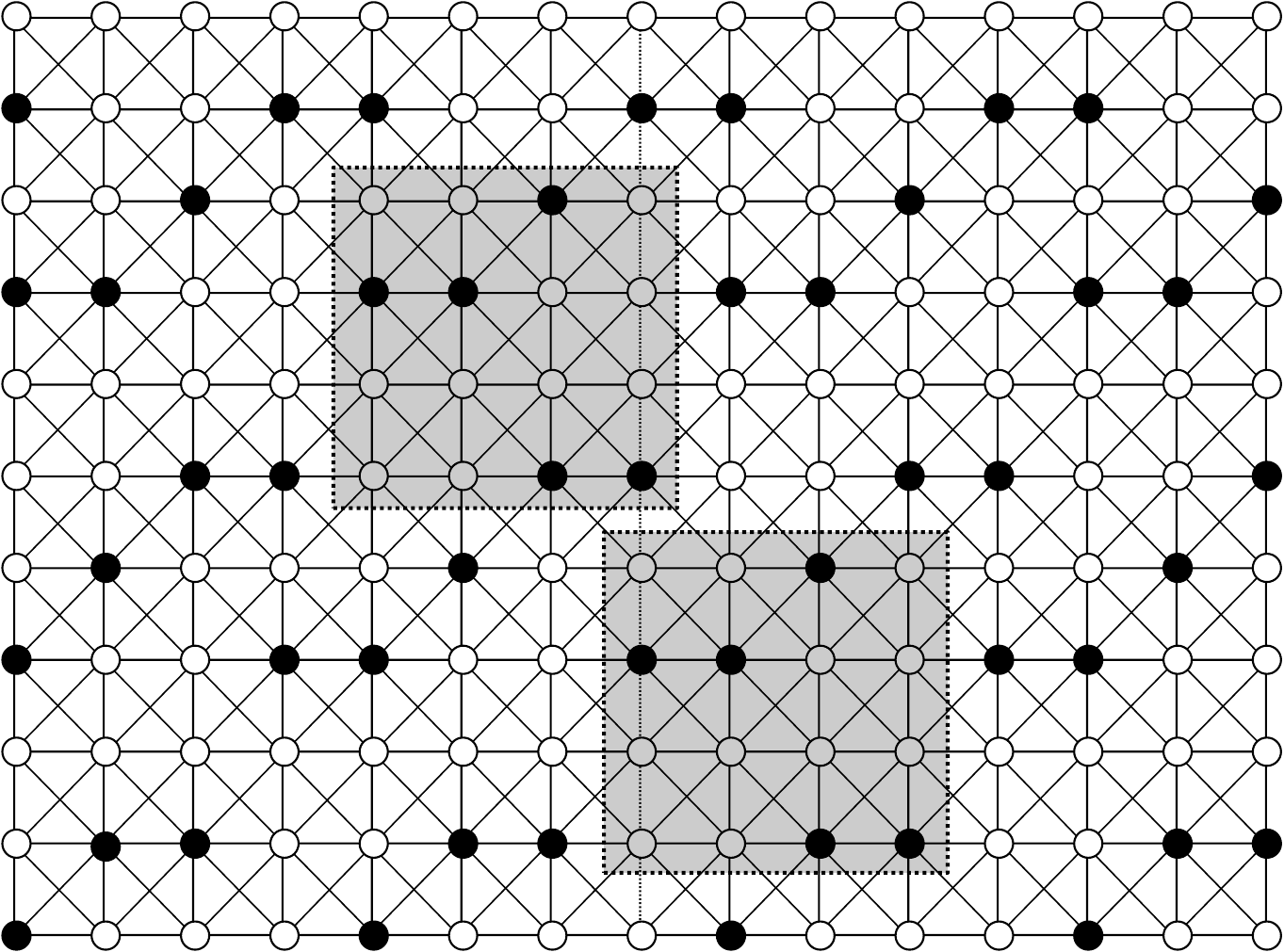} \\ (a) & (b)
        \\ \\
        \includegraphics[width=0.375\textwidth]{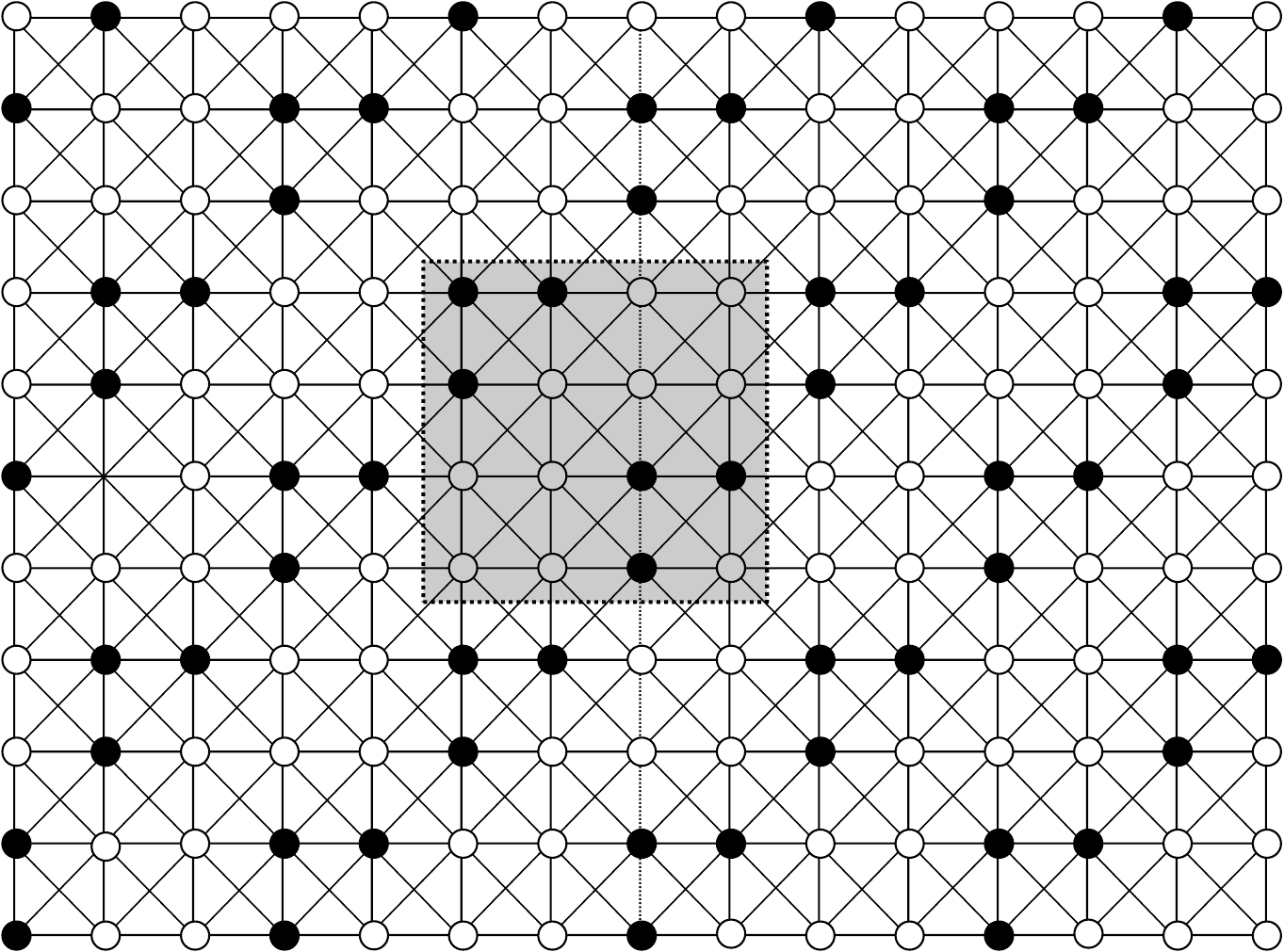} &
        \includegraphics[width=0.375\textwidth]{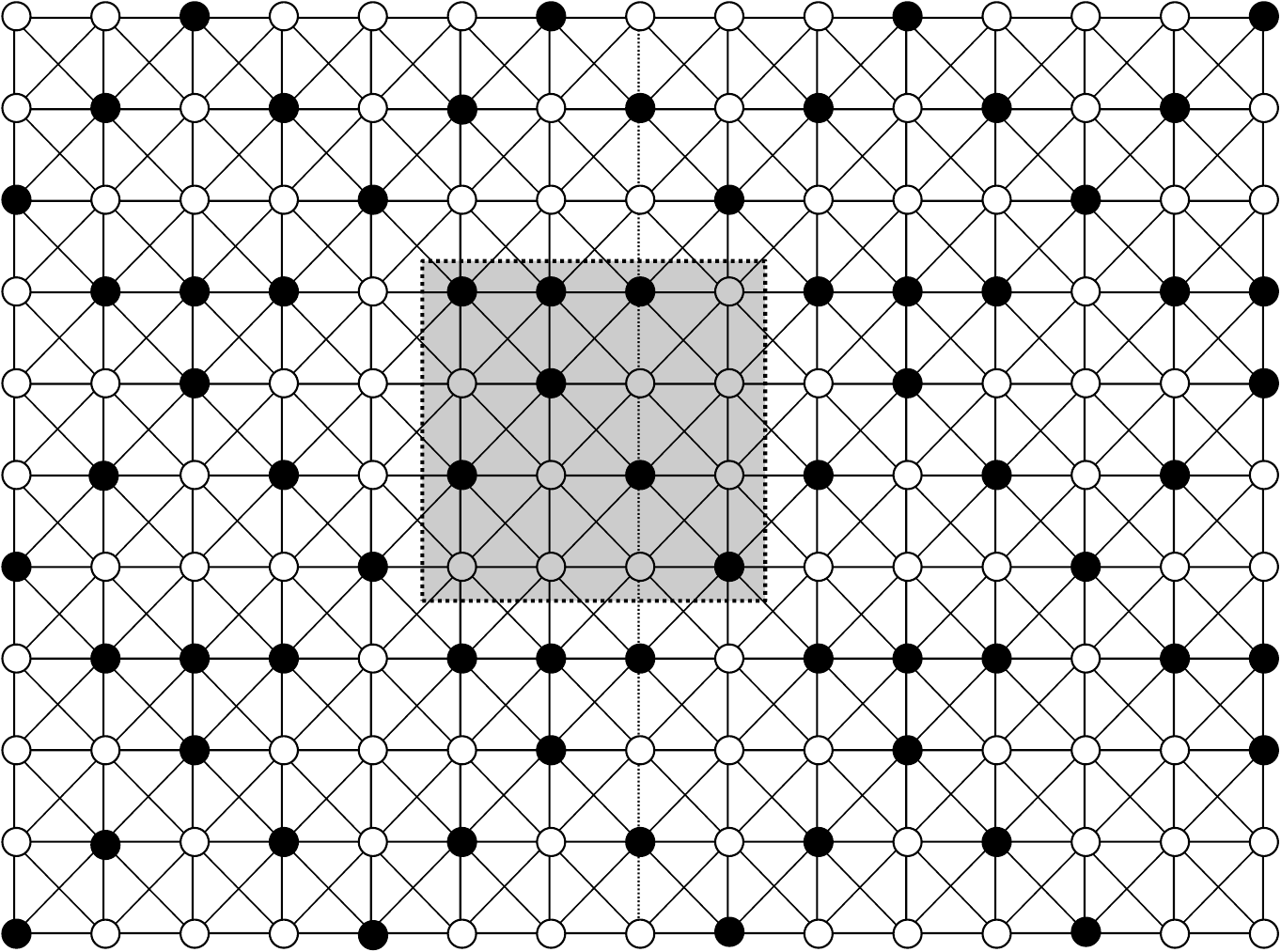} \\ (c) & (d)
    \end{tabular}
    \caption{Our best constructions of LD (a), RED:LD (b), DET:LD (c), and ERR:LD (d) sets on the King's grid. Shaded vertices denote detectors.}
    \label{fig:king-solns}
\end{figure}

In this section we will determine upper and lower bounds for the optimal (minimum) densities of RED:LD, DET:LD, and ERR:LD sets on the infinite king grid, $\textrm{K}$.
Firstly, for the upper bound, we have constructed LD, RED:LD, DET:LD, and ERR:LD sets on the infinite king grid with densities $\frac{1}{5}$, $\frac{5}{16}$, $\frac{6}{16}$, and $\frac{7}{16}$, respectively, as shown in Figure~\ref{fig:king-solns} (note that Figure~\ref{fig:king-solns}~(a) was found to be isomorphic to a previous solution by Honkala and Laihonen \cite{honk05c}).
One may verify that each of these constructions are valid detection systems using Theorems \ref{theo:ld-char}, \ref{theo:red-ld-char}, \ref{theo:det-ld-char}, and \ref{theo:err-ld-char}, respectively.

\begin{figure}[ht]
    \centering
    \includegraphics[width=0.3\textwidth]{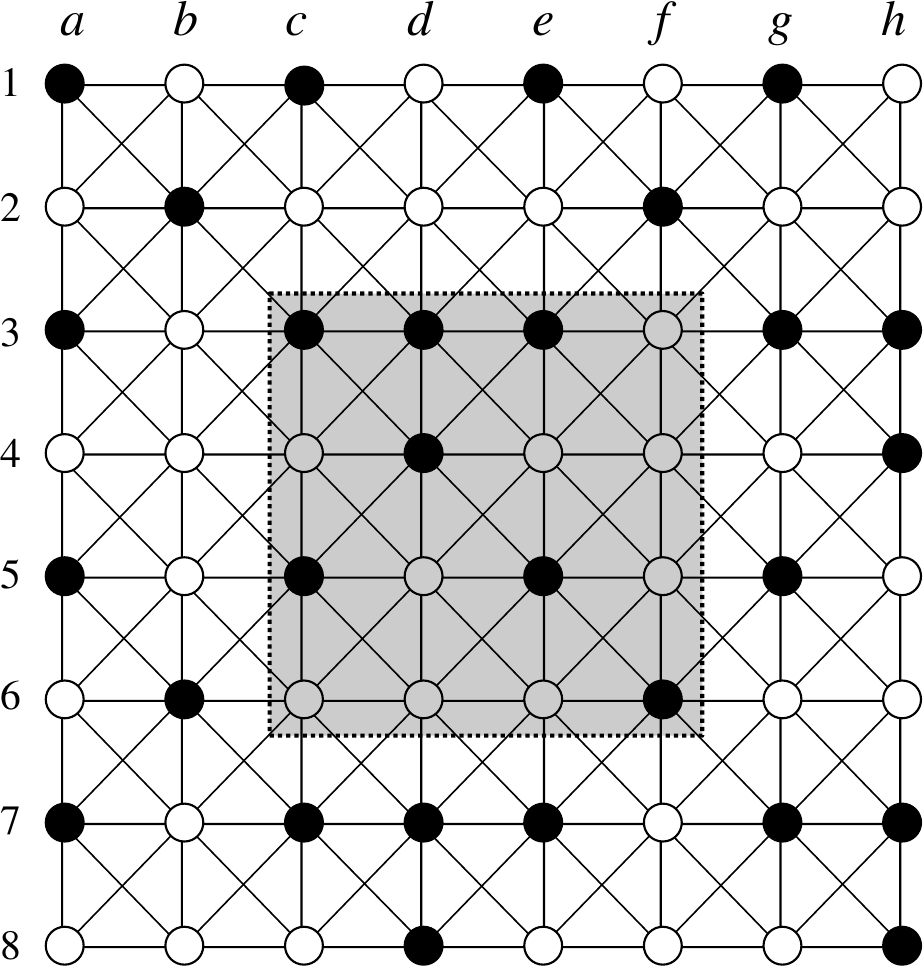}
    \caption{ERR:LD set to demonstrate Theorem~\ref{theo:err-ld-char}}
    \label{fig:err-ld-one-tile}
\end{figure}

As an example, we will show Figure~\ref{fig:king-solns}~(d) is an ERR:LD set; we will use the labeling from Figure~\ref{fig:err-ld-one-tile}.
We observe that every vertex in $V(\textrm{K})$ is 3-dominated, so $S$ satisfies property i) of Theorem~\ref{theo:err-ld-char}.
Now, we will check every vertex pair $u,v \in V(G)$ with $u \neq v$ is 3-distinguished.
If $d(u,v) \ge 3$, then clearly the pair is 6-distinguished because no detector can dominate both of them, and so the vertex pair $(u,v)$ is distinguished; thus, we need only consider when $d(u,v) \le 2$.
Because every local (shaded) tile will be identical to any other, we can assume that $u$ is in the tile.
We then consider $v$ being any vertex such that $d(u,v) \le 2$.
If every pair is distinguished, then we will have satisfied all four properties of Theorem~\ref{theo:err-ld-char}, meaning $S$ is an ERR:LD set.
For example, take the detector pair $(v_{c3}, v_{d3})$: aside from themselves, they are dominated by $v_{b2}$, $v_{e3}$, and $v_{d4}$; however $v_{d4}$ is common to both.
Thus, $v_{c3}$ and $v_{c4}$ are 2-distinguished; they need only be 1-distinguished, so they are distinguished.
As another example, take the detector pair $(v_{e5},v_{d7})$: they have no common detectors, and are clearly 5-distinguished, so they are distinguished.
Consider the non-detector pair $(v_{f4},v_{f5})$; we see they are 3-distinguished by $v_{e3}$, $v_{g3}$, and $v_{f6}$, and therefore meet the minimum requirements to be distinguished.
We also see the non-detector pair $(v_{f5},v_{g6})$ are 3-distinguished by $v_{e5}$, $v_{g7}$, and $v_{h7}$, and are thus distinguished.
Consider the mixed pair $(v_{c5},v_{d5})$; we find they are 2-distinguished by $v_{b6}$ and $v_{e5}$, and thus satisfy the minimum requirements to be distinguished.
As another example, we see the mixed pair $(v_{d2},v_{d4})$ is 4-distinguished by $v_{c1}$, $v_{e1}$, $v_{c5}$, and $v_{e5}$, and thus is distinguished.
It can be shown that all pairs are distinguished, so $S$ is an ERR:LD set.
Thus we have an upper bound for the optimal density: $\textrm{ERR:LD}(\textrm{K}) \le \frac{7}{16}$.

Similarly, we can verify that the tilings in (a), (b), and (c) from Figure~\ref{fig:king-solns} satisfy the requirements of Theorems \ref{theo:ld-char}, \ref{theo:red-ld-char}, and \ref{theo:det-ld-char}, so we have upper bounds for the optimal densities: $\textrm{RED:LD}(\textrm{K}) \le \frac{5}{16}$ and $\textrm{DET:LD}(\textrm{K}) \le \frac{6}{16}$.
Figure~\ref{fig:king-solns} shows an upper bound of $\frac{1}{5}$ for $\textrm{LD}(\textrm{K})$, which has already been proven to be optimal \cite{honk05c}; it is included for comparison purposes.

With upper bounds for each parameter now established, we will proceed to explore lower bounds for these sets.
If $S$ is a RED:LD set, then Theorem~\ref{theo:red-ld-char} guarantees every $v \in V(\textrm{K})$ must be 2-dominated, so $sh[v] \le \frac{1}{2}$; thus, $\forall x \in S$, $sh(x) \le 9 \times \frac{1}{2} = \frac{9}{2}$, giving us a simple lower bound: $\textrm{RED:LD}(\textrm{K}) \ge \frac{2}{9}$.
Similarly, by Theorems \ref{theo:det-ld-char} and \ref{theo:err-ld-char}, we find that $\textrm{DET:LD}(\textrm{K}) \ge \frac{2}{9}$ and $\textrm{ERR:LD}(\textrm{K}) \ge \frac{3}{9} = \frac{1}{3}$.
However, these simple lower bounds can be greatly improved by applying the discharging technique we have described.
In what follows, we will use this technique to improve these bounds for RED:LD and ERR:LD to $\frac{3}{11}$ and $\frac{15}{38}$, respectively.


\newpage
\begin{theorem}\label{theo:sh-err-ld-king}
Let $S \subseteq V(G)$ be an ERR:LD set for the king's grid.
Then, the average share of all vertices in $S$ is at most $\frac{38}{15}$.
\end{theorem}
\cbeginproof

\begin{figure}[ht]
    \centering
    \begin{tabular}{c@{\hspace{4em}}c@{\hspace{4em}}c@{\hspace{4em}}c}
        \includegraphics[width=0.125\textwidth]{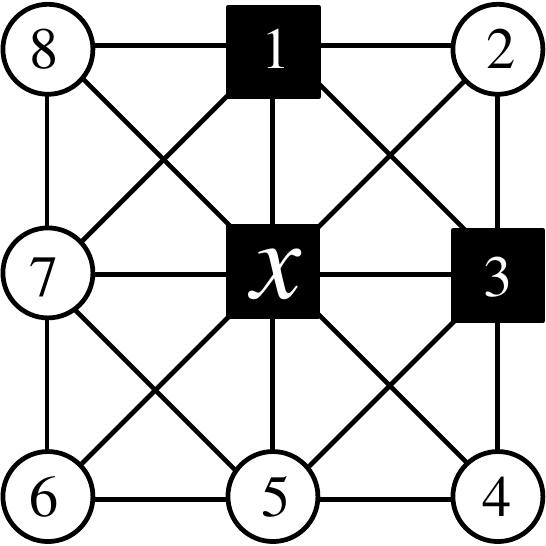} & \includegraphics[width=0.125\textwidth]{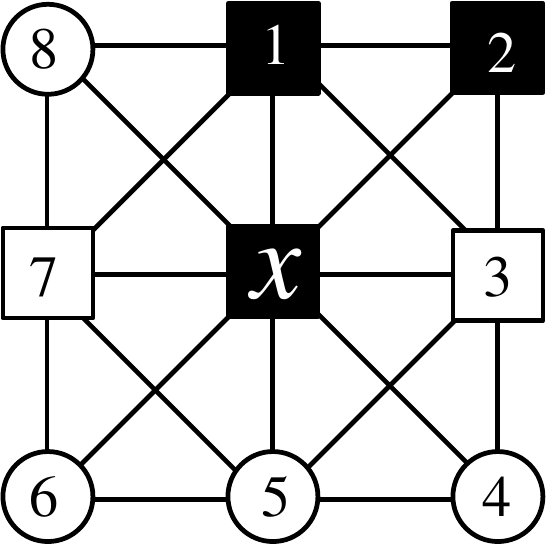} & \includegraphics[width=0.125\textwidth]{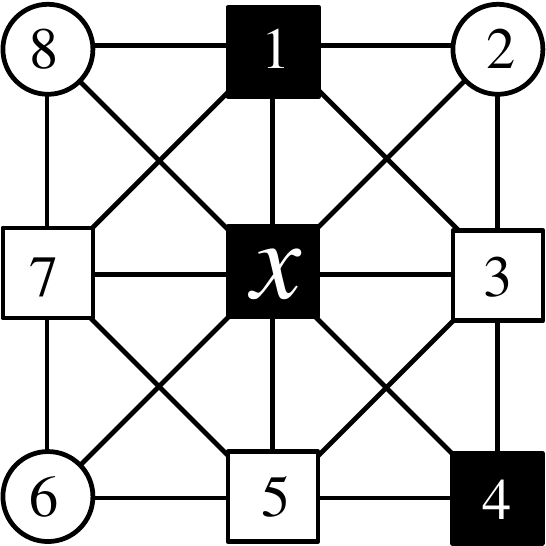} & \includegraphics[width=0.125\textwidth]{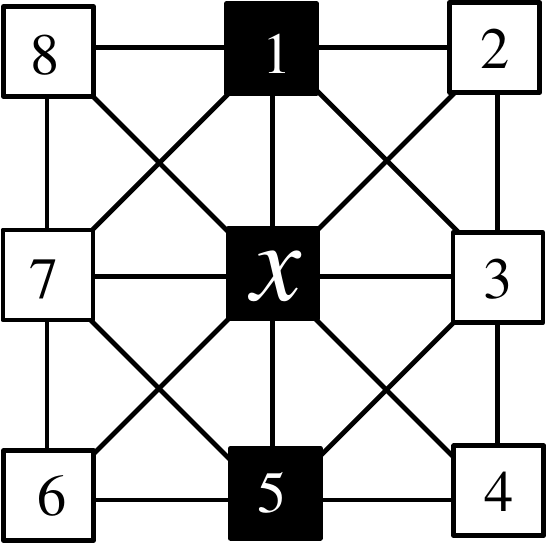} \\ \\
        
        \includegraphics[width=0.125\textwidth]{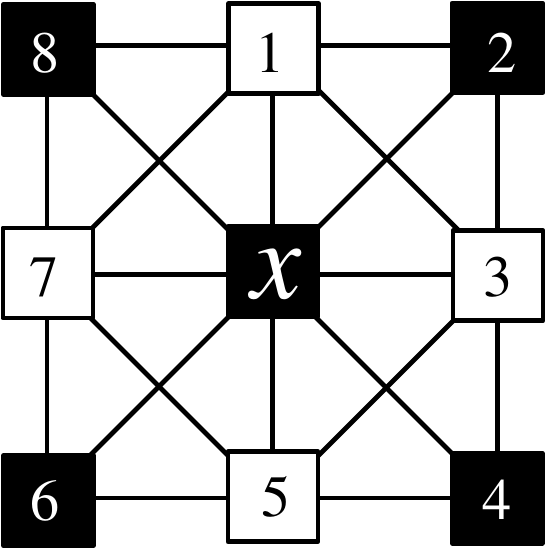} & \includegraphics[width=0.125\textwidth]{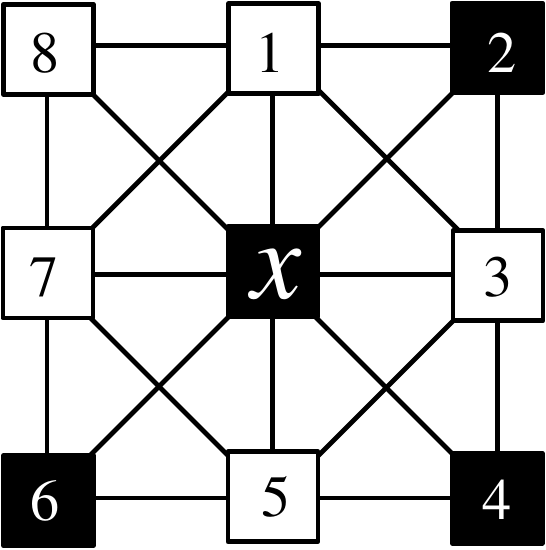} & \includegraphics[width=0.125\textwidth]{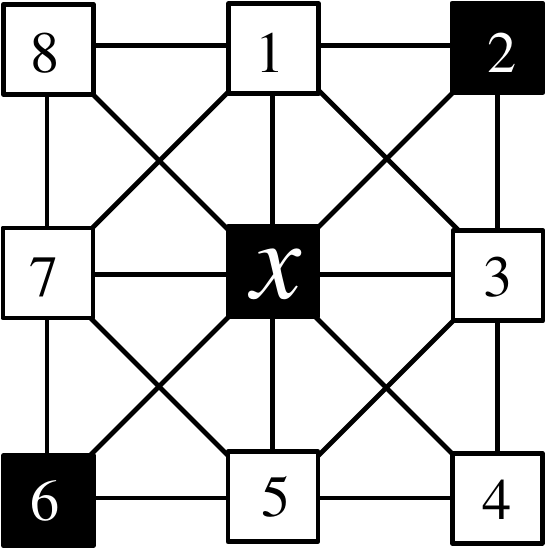} & \includegraphics[width=0.125\textwidth]{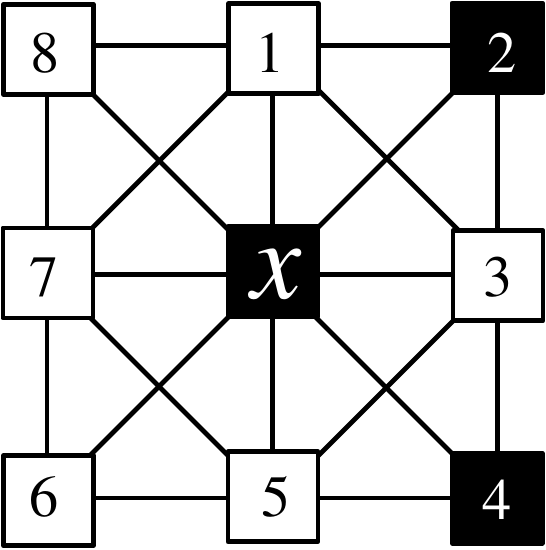} \\ 
    \end{tabular}
    \caption{Configurations around a detector vertex $x$.}
    \label{fig:err-ld-cases}
\end{figure}

Let $S \subseteq V(G)$ be an ERR:LD set on $\textrm{K}$, let the notations $sh$, $dom$, and $k$-dominated implicitly use $S$, and let a vertex pair being \emph{distinguished} denote satisfying the relevant pair property of Theorem~\ref{theo:err-ld-char}.
Let $x \in S$ be a detector vertex; Theorem~\ref{theo:err-ld-char} yields that $|N[x] \cap S| \ge 3$.
There are eight non-isomorphic ways to at least 3-dominate $x$, as shown in Figure~\ref{fig:err-ld-cases}.
Square vertices represent detectors (shaded) or non-detectors (non-shaded), and round vertices denote unknowns.
Note that vertices are labeled with numbers such as $k$, but will be referenced in the text as $v_k$.

\begin{wrapfigure}{r}{0.23\textwidth}
    \centering
    \vspace{-2em}
    \includegraphics[width=0.2\textwidth]{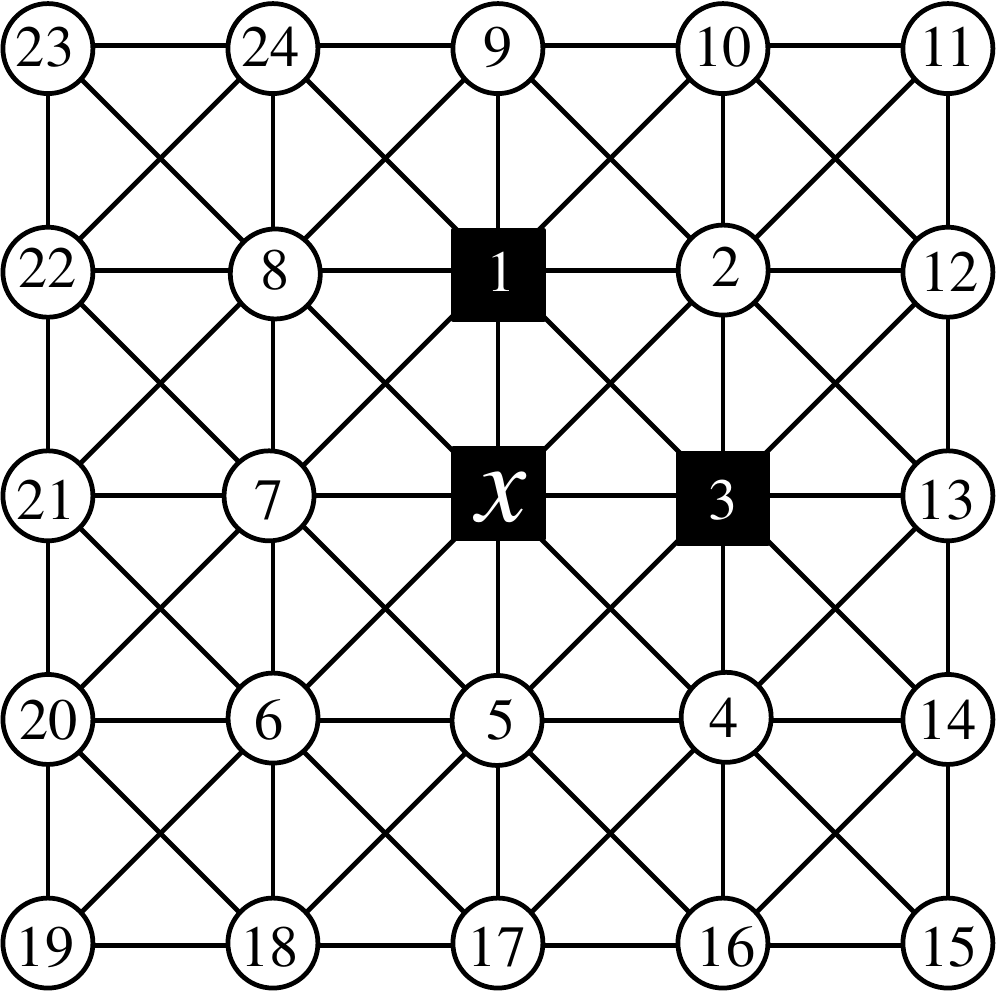}
    \caption{L shape}
    \label{fig:err-ld-case-elshape}
\end{wrapfigure}
\newcase{L shape}
First, we consider the partial shares of vertices $A = \{x,v_1,v_2,v_3\}$.
We require $|A \cap D_{4+}| \ge 2$ for the vertices in $A$ to be distinguished.
Suppose $v_2$ is a non-detector.
Let $d = \min_{w \in \{x,v_1,v_3\}}{dom(w)}$.
If $d = 3$ then $dom(v_2) \ge 5$ to be distinguished; if $d = 4$ then $dom(v_2) \ge 4$ to be distinguished; otherwise $d \ge 5$ and $dom(v_2) \ge 3$.
Thus, $sh[xv_1v_2v_3] \le \max \{ \sigma_{3445}, \sigma_{3445}, \sigma_{4444}, \sigma_{5553}\} = \frac{31}{30}$.
Otherwise, $v_2$ is a detector; in order to distinguish it from the other three detectors, $x$, $v_1$, and $v_3$, we need $sh[xv_1v_2v_3] \le \sigma_{4555} < \frac{31}{30}$.
In either case, $sh[xv_1v_2v_3] \le \frac{31}{30}$.

Next, we consider the partial shares of $v_4$ and $v_5$.
If both are detectors, then $sh[v_4v_5] \le \sigma_{45}$ to be distinguished.
If both are non-detectors, then $sh[v_4v_5] \le \sigma_{34}$ to be distinguished.
Otherwise, we let $v_1$ be unknown to allow symmetry between $v_4$ and $v5$; then by symmetry let $\{v_4,v_5\} \cap S = \{v_4\}$; then $sh[v_4v_5] \le \max \{\sigma_{35},\sigma_{44}\}$.
Thus, $sh[v_4v_5] \le \sigma_{34}$, and by symmetry, $sh[v_7v_8] \le \sigma_{34}$ as well.
Therefore, $sh(x) = sh[xv_1v_2v_3] + sh[v_4v_5] + sh[v_6] + sh[v_7v_8] \le \sigma_{3445} + \sigma_{34}  + \sigma_3 + \sigma_{34} = \frac{38}{15}$, and we are done.

\FloatBarrier
\begin{wrapfigure}[6]{r}{0.2\textwidth}
    \centering
    \vspace{-2em}
    \includegraphics[width=0.2\textwidth]{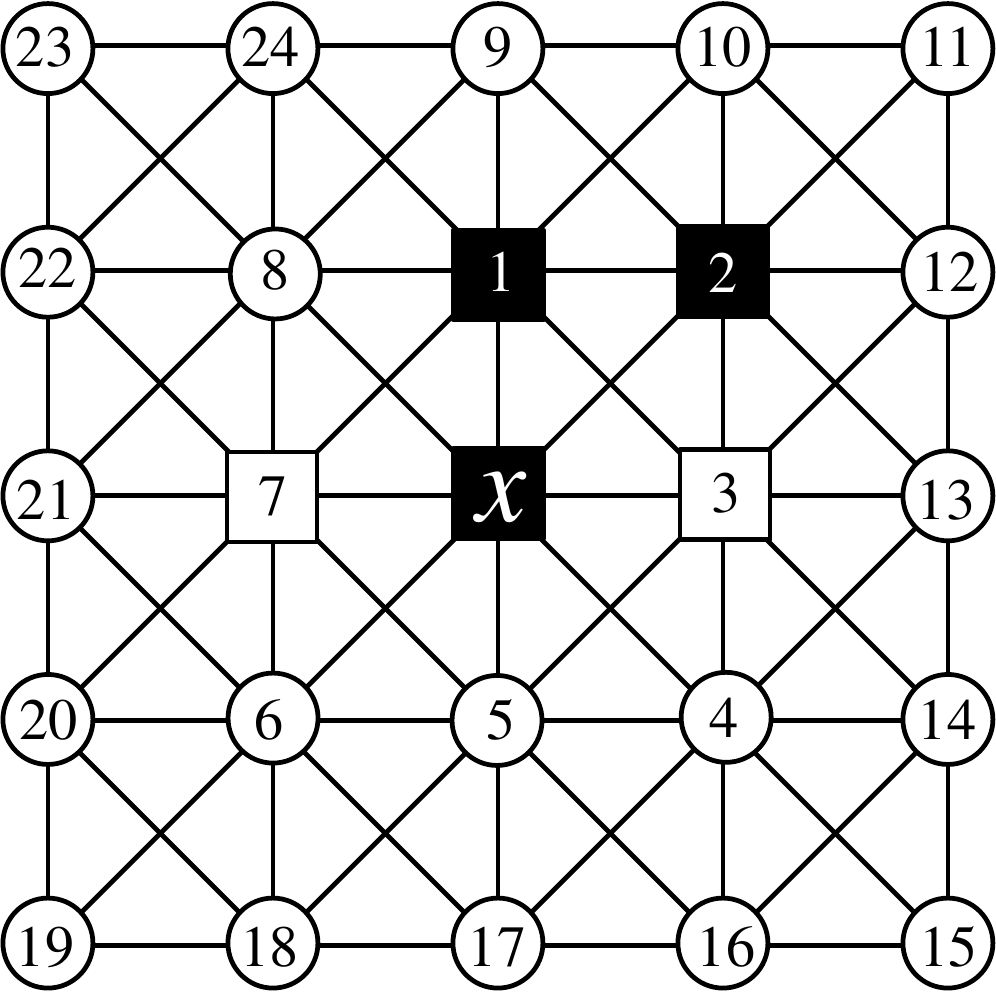}
    \caption{Triangle}
    \label{fig:err-ld-case-triangle}
\end{wrapfigure}

\newcase{Triangle}
Next, we consider the case where $x$ is dominated by two neighbors in a ``triangle'' shape, as in Figure~\ref{fig:err-ld-case-triangle}.
We will assume that $\{v_3,v_7\} \cap S = \varnothing$, as otherwise we fall into the previous case.

We consider $sh[xv_1v_2v_3]$, $sh[v_7v_8]$, and $sh[v_4v_5v_6]$.
First, if $dom(v_3) = 3$, then $\{x,v_1,v_2\} \subseteq D_{5+}$ to distinguish from $v_3$, so $sh[xv_1v_2v_3] \le \sigma_{3555}$; if $dom(v_3) = 4$, then $sh[xv_1v_2v_3] \le \sigma_{4444}$ to be distinguished; otherwise $dom(v_3) \ge 5$, in which case $sh[xv_1v_2v_3] \le \sigma_{3445}$.
Therefore, $sh[xv_1v_2v_3] \le \sigma_{3445}$. 
If $v_8 \in S$, then $sh[v_7v_8] \le \max \{\sigma_{35}, \sigma_{44}\} < \sigma_{34}$; otherwise, $sh[v_7v_8] \le \sigma_{34}$.
If $\{v_4, v_5, v_6\} \cap D_{4+} \neq \varnothing$, then $sh(x) = sh[xv_1v_2v_3] + sh[v_4v_5v_6] + sh[v_7v_8] \le \sigma_{3445} + \sigma_{334} + \sigma_{34} = \frac{38}{15}$ and we will be done.

Otherwise, we assume $\{v_4, v_5, v_6\} \subseteq D_3$; we consider the following sub-cases depending on whether these three vertices are detectors or not.
\textbf{Subcase~1:} If $\{v_4,v_5,v_6\} \cap S = \varnothing$, then $v_5$ needs two additional detectors, but at least one of them will be in $N(v_4)$; thus, $v_4$ and $v_5$ are not distinguished, a contradiction.
\textbf{Subcase~2:} If $\{v_4,v_5,v_6\} \cap S = \{v_5\}$, then $v_4$ and $v_6$ are not distinguished, a contradiction.
\textbf{Subcase~3:} If $\{v_4,v_5,v_6\} \cap S = \{v_6\}$, then $v_5$ needs one additional detector, $w \in S$; if $w \in N(v_4)$ then $\{v_4,v_5\}$ are not distinguished, a contradiction; otherwise $w = v_{18}$ and $\{v_5,v_6\}$ are not distinguished, a contradiction.
Similarly, if $\{v_4,v_5,v_6\} \cap S = \{v_4\}$ then we reach a contradiction.
\textbf{Subcase~4:} If $\{v_4,v_5,v_6\} \cap S = \{v_4,v_5\}$ then $v_4$ and $v_5$ are not distinguished, a contradiction; similarly, if $\{v_4,v_5,v_6\} \cap S = \{v_5,v_6\}$ then $v_5$ and $v_6$ are not distinguished, contradiction.
\textbf{Subcase~5:} If $\{v_4,v_5,v_6\} \cap S = \{v_4,v_6\}$, consider the third dominator for $v_6$.
If $v_{17} \in S$ or $v_{18} \in S$ is the third dominator, we get $dom(v_5) \ge 4$, a contradiction.
If $v_{20} \in S$ or $v_{21} \in S$ is the third dominator, then we need $v_8 \in S$ or $v_{22} \in S$ to distinguish $v_7$ from $v_6$, which in turn will make $dom(v_7) \ge 5$ and $dom(v_8) \ge 4$; thus $sh(x) = sh[xv_1v_2v_3] + sh[v_4v_5v_6] + sh[v_7v_8] \le \sigma_{3445} + \sigma_{333} + \sigma_{45} < \frac{38}{15}$.
Otherwise, $v_{19} \in S$ is the third dominator.
Then we need $v_8 \in S$ or $v_{22} \in S$ to distinguish $v_7$ from $v_5$, which will make $dom(v_7) \ge 4$ and $dom(v_8) \ge 4$; thus $sh(x) = sh[xv_1v_2v_3] + sh[v_4v_5v_6] + sh[v_7v_8] \le \sigma_{3445} + \sigma_{333} + \sigma_{44} = \frac{38}{15}$, completing the proof of the triangle case.

\begin{wrapfigure}{r}{0.2\textwidth}
    \centering
    \vspace{-1em}
    \includegraphics[width=0.2\textwidth]{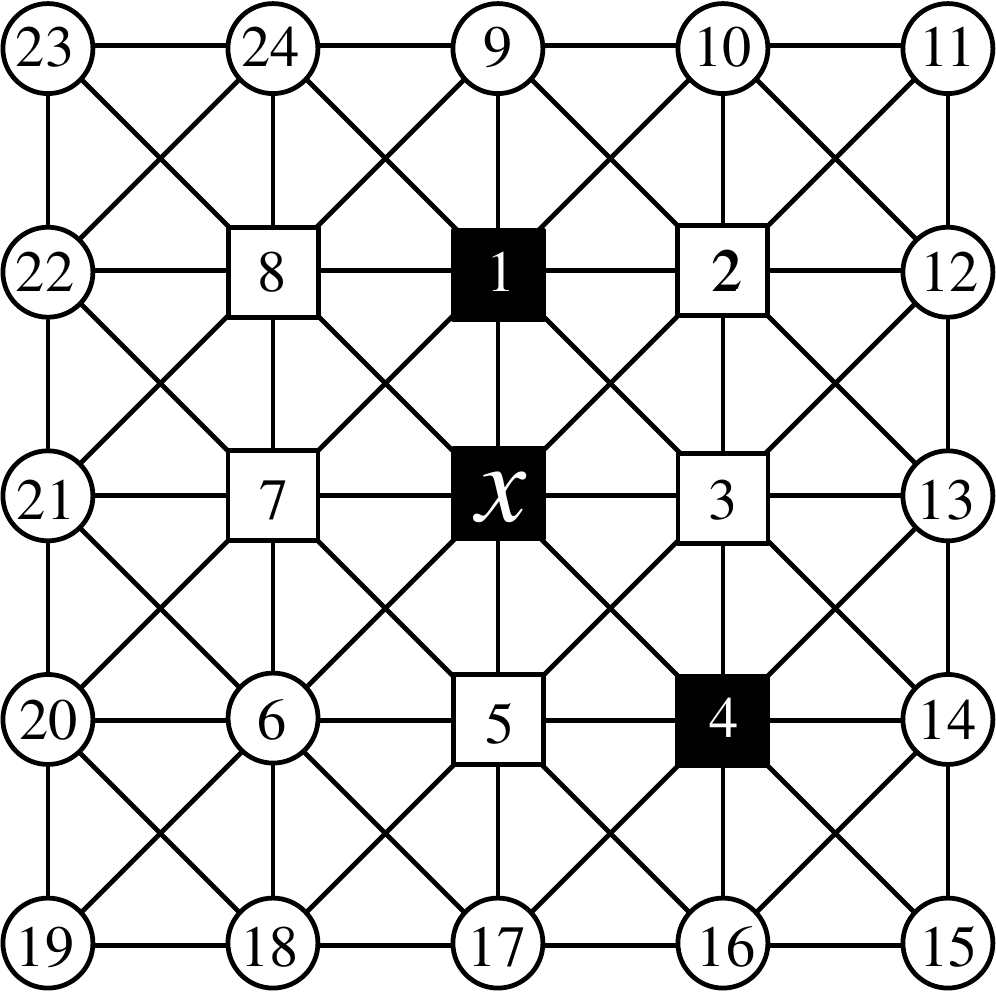}
    \caption{Bent}
    \label{fig:err-ld-case-bent}
\end{wrapfigure}

\newcase{``Bent''}
Consider the case where $x$ is dominated by two neighbors in a ``bent'' shape, as in Figure \ref{fig:err-ld-case-bent}.
We will assume that $\{v_2, v_3, v_5, v_7, v_8\} \cap S = \varnothing$, as otherwise we fall into previous cases.

First, we consider when $v_6 \in S$, then $dom(x) = 4$.
We observe that $\{v_3,v_5,v_7\} \subseteq D_{4+}$ to be distinguished from $x$.
We now consider how to 4-dominate $v_7$.
Suppose $v_{20} \notin S$, then $\{v_{21},v_{22}\} \cap S \neq \varnothing$, and we observe that $dom(v_8) \ge 5$ to be distinguished from $v_7$, so $sh(x) \le \frac{38}{15}$ and we are done.
Thus, we assume that $v_{20} \in S$.
If $\{v_{21},v_{22}\} \cap S \neq \varnothing$ then again $dom(v_8) \ge 5$ to be distinguished from $v_7$ and we are done; so we assume $\{v_{21},v_{22}\} \cap S = \varnothing$.
By symmetry, we can assume that $\{v_{12},v_{13},v_{14}\} \cap S = \{v_{14}\}$.
If $\{v_2,v_8\} \subseteq D_{4+}$ we are done, so without loss of generality we assume $dom(v_8) = 3$.
If $v_{23} \in S$ then $v_1$ and $v_2$ cannot be distinguished, a contradiction.
Otherwise $v_{24} \in S$ or $v_9 \in S$; in either case, $v_1$ and $v_8$ cannot be distinguished, a contradiction.
Therefore, $sh(x) \le \frac{38}{15}$.

Finally, we consider when $v_6 \notin S$.
We observe that $dom(v_3) \ge 5$ to be distinguished from $x$, which in turn requires $dom(v_2) \ge 4$ to distinguish it from $v_3$; thus, $sh[xv_2v_3] \le \sigma_{345}$.
Next, we consider $sh[v_1v_7v_8]$ based on $dom(v_7)$.
If $v_7$ is 5-dominated, then we need $dom(v_8) \ge 6$ to distinguish it from $v_7$.
If $v_7$ is 4-dominated with $\{v_{20},v_{21}\}$ or $\{v_{20},v_{22}\}$, we need $dom(v_8) \ge 5$ to distinguish from $v_7$; otherwise, $v_7$ is 4-dominated with $\{v_{21},v_{22}\}$ and we need $dom(v_8) \ge 7$ to distinguish from $v_7$.
If $v_7$ is 3-dominated with $v_{21}$ or $v_{22}$, we need $dom(v_8) \ge 6$ to distinguish from $v_7$; otherwise, $v_7$ is 3-dominated with $v_{21}$ and we need $dom(v_8) \ge 4$, which in turn requires $dom(v_1) \ge 4$.
From the three cases of $dom(v_7)$, we have $sh[v_1v_7v_8]\le \sigma_{344}$.
If $\{v_4,v_5,v_6\} \cap D_{4+} \neq \varnothing$ then we will be done because $sh(x) = sh[xv_2v_3] +  sh[v_4v_5v_6] + sh[v_1v_7v_8] \le \sigma_{345} + \sigma_{334} + \sigma_{344} \le \frac{38}{15}$.
Otherwise, we assume $\{v_4,v_5,v_6\} \subseteq D_3$.
Then, we observe $\{v_{15},v_{16},v_{17}\} \cap S = \varnothing$ because $v_4$ is already 3-dominated, which in turn requires $v_{18} \in S$ to 3-dominate $v_5$.
We find that $v_5$ and $v_6$ are not distinguished, a contradiction, completing the proof for Case ``bent.''

\begin{wrapfigure}[8]{r}{0.2\textwidth}
    \centering
    \vspace{-2.5em}
    \includegraphics[width=0.2\textwidth]{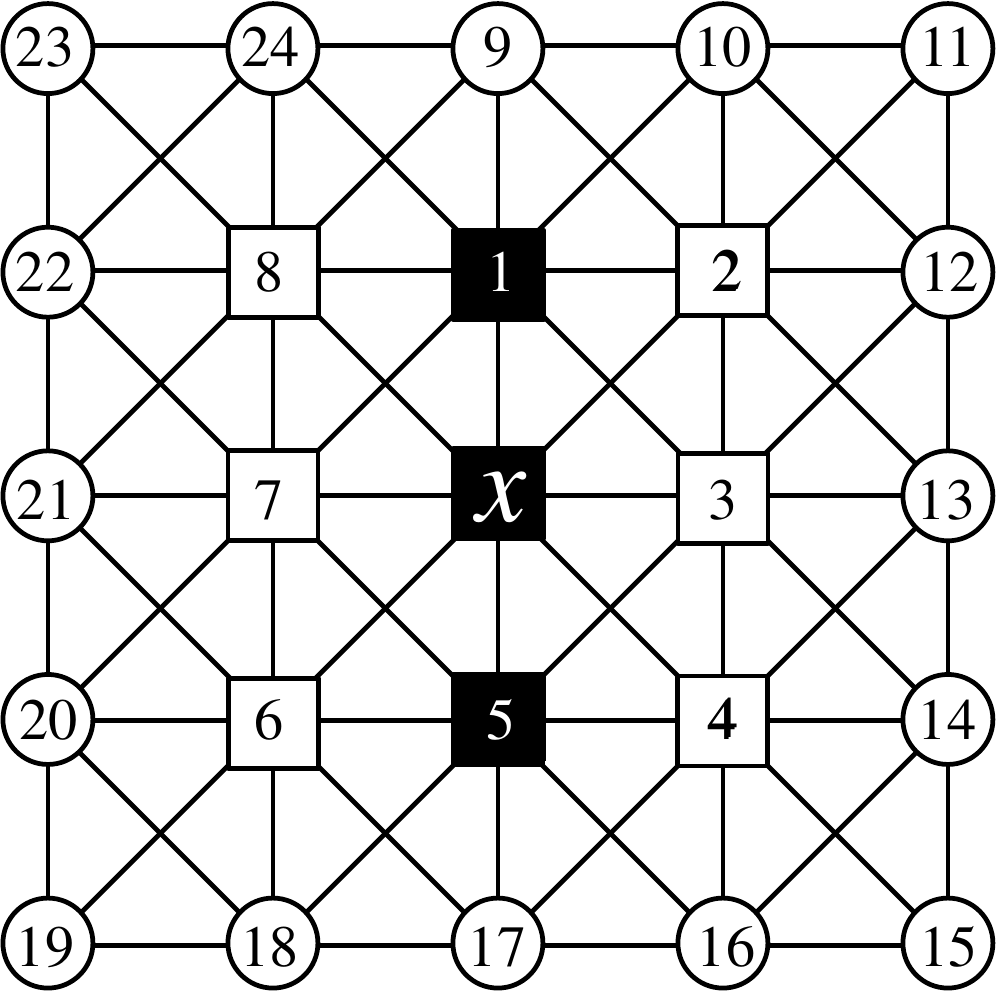}
    \caption{Vertical}
    \label{fig:err-ld-case-vertical}
\end{wrapfigure}

\newcase{Vertical}
Consider the case where $x$ is dominated by two neighbors in a ``vertical'' shape, as in Figure \ref{fig:err-ld-case-vertical}.
We can assume that $\{v_2, v_3, v_4, v_6, v_7, v_8\} \cap S = \varnothing$, as otherwise we fall into previous cases.

We observe that $dom(v_3) \ge 5$ to be distinguished from $x$.
Further, $\{v_2,v_4\} \subseteq D_{4+}$ to be distinguished from $v_3$, thus $sh[v_2v_3v_4] \le \sigma_{454}$.
By symmetry, we have $sh[v_6v_7v_8]\le \sigma_{454}$, as well.
Therefore, $sh(x) = sh[xv_1v_5]+ sh[v_2v_3v_4] +sh[v_7v_8v_9]\le \sigma_{333} + \sigma_{454} + \sigma_{454} \le \frac{38}{15}$, and we are done.

\begin{wrapfigure}{r}{0.2\textwidth}
    \centering
    \vspace{-1em}
    \includegraphics[width=0.2\textwidth]{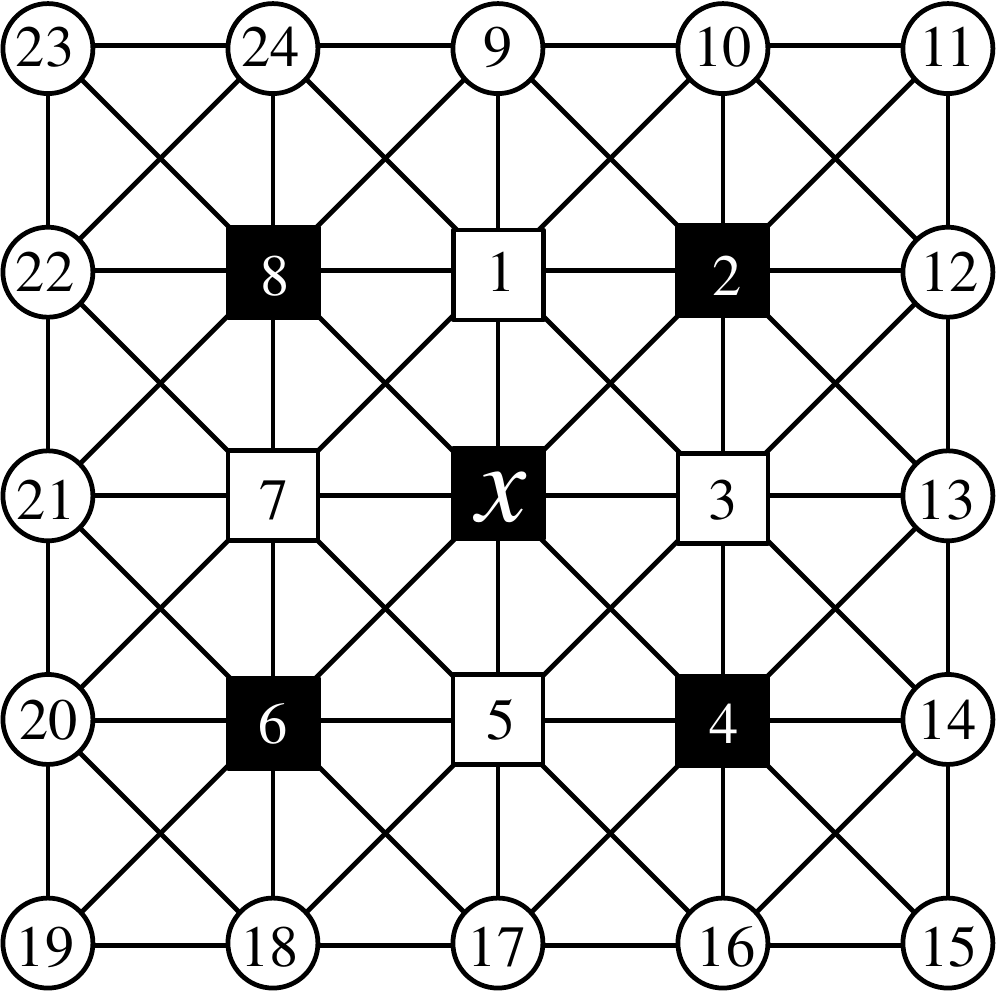}
    \caption{X-shape}
    \label{fig:err-ld-case-xshape}
\end{wrapfigure}

\newcase{``X'' Shape}
Consider the case where $x$ is dominated by four neighbors in a ``X'' shape, as in Figure \ref{fig:err-ld-case-xshape}.
We can assume that $\{v_1,v_3,v_5,v_7\} \cap S = \varnothing$, as otherwise we fall into previous cases.

Suppose $\{v_1,v_3,v_5,v_7\} \subseteq D_{4+}$, then $sh(x) \le \frac{38}{15}$ and we are done; otherwise, without loss of generality we assume $dom(v_7) = 3$.
Then $\{v_1,v_5\} \subseteq D_{4+}$ to be distinguished from $v_7$.
Suppose $dom(v_3) \ge 4$.
If $v_{14} \notin S$, then $dom(v_2) \ge 4$ to be distinguished from $v_3$ and we are done; thus, we assume with symmetry that $\{v_{12},v_{14}\} \subseteq S$.
If $v_{13} \in S$, then $dom(v_2) \ge 4$ and we are done, so we assume $v_{13} \notin S$.
If $dom(v_2) \ge 4$ we are done, so we assume $\{v_9,v_{10},v_{11}\} \cap S = \varnothing$, requiring $v_{24} \in S$ to 3-dominate $v_1$.
If $dom(v_8) \ge 4$ we are done, otherwise we find $v_1$ and $v_8$ are not distinguished, a contradiction.
Now, we can assume $\{v_3,v_7\} \subseteq D_3$ and examine $dom(v_1)$ and $dom(v_5)$.
Without loss of generality we can assume $dom(v_1) \le dom(v_5)$; then $(dom(v_1), dom(v_5)) \in \{(4, 4), (4, 5), (4, 6), (5, 5), (5, 6), (6, 6)\}$ and consider the following four sub-cases: where $dom(v_5) = 6$, and the other three possibilities.

First, if $dom(v_5) = 6$, then it requires $\{v_{16},v_{17},v_{18}\} \subseteq S$, which implies that $\{v_4,v_6\} \subseteq D_4$, so $sh(x) = sh[xv_2v_3v_7v_8] + sh[v_1v_5] + sh[v_4v_6]  \le \sigma_{53333}+ \sigma_{46} +\sigma_{44} \le \frac{38}{15}$.
Second, assume $dom(v_1) = dom(v_5) = 4$.
If $v_9 \in S$ to 4-dominate $v_1$, then $\{v_2,v_8\} \subseteq D_{4+}$ to distinguish from $v_1$.
Otherwise, $v_{24} \in S$ or $v_{10} \in S$, then $\{v_2,v_8\} \cap D_{4+} \neq \varnothing$ to distinguish from $v_1$.
Hence, $sh[v_1v_2v_8] \le \sigma_{344}$ and by symmetry, we also have $sh[v_4v_5v_6] \le \sigma_{344}$.
Therefore, $sh(x) = sh[xv_3v_7] + sh[v_1v_2v_8] + sh[v_4v_5v_6] \le \sigma_{533} +\sigma_{344} + \sigma_{344} = \frac{38}{15}$ and we are done.

Next, assume $dom(v_1)=4$ and $dom(v_5)=5$.
If $v_9 \in S$, then we need $\{v_2,v_8\} \subseteq D_{4+}$ for $v_1$ to be distinguished from $v_2$ and $v_8$, and we get $sh(x) = sh[xv_3v_7] + sh[v_1v_2v_8] + sh[v_4v_5v_6] \le \sigma_{533} +\sigma_{444} + \sigma_{353} = \frac{149}{60} < \frac{38}{15}$, so we are done; thus, we assume $v_9 \notin S$.
We need  $v_{10} \in S$ or $v_{24} \in S$ to 4-dominate $v_1$; by symmetry, let $v_{24} \in S$, which requires $v_{11} \in S$ to 3-dominate $v_2$ and $v_{23} \in S$ to distinguish $v_8$ from $v_1$.
If $v_{17} \in S$, then we need $v_{16} \in S$ or $v_{18} \in S$ to 5-dominate $v_5$; by symmetry, let $v_{16} \in S$, then $dom(v_4) \ge 5$ to be distinguished from $v_5$ and we are done.
Thus, we assume $v_{17} \notin S$, requiring $\{v_{16}, v_{18}\} \subseteq S$ to 5-dominate $v_5$.
If $v_{15} \in S$ or $v_{19} \in S$, then $sh(x) = sh[xv_3v_7] + sh[v_1v_2v_8] + sh[v_4v_5v_6] \le \sigma_{533} +\sigma_{434} + \sigma_{453} \le \frac{38}{15}$, and we are done. 
Therefore, we assume $v_{15} \notin S$ and $v_{19} \notin S$, which leads us to the first problem configuration, which will be handled separately.

Lastly, assume $dom(v_1)=5$ and $dom(v_5)=5$.
Similarly to the previous case, if $v_{17} \in S$, then $sh(x) \le \frac{38}{15}$ and we are done, so with symmetry we assume $v_{17} \notin S$ and $v_9 \notin S$.
Therefore, we require $\{v_{10}, v_{24}, v_{16}, v_{18}\} \subseteq S$ to 5-dominate $v_1$ and $v_5$.
If $\{v_{11}, v_{15}, v_{19}, v_{23}\} \cap S \neq \varnothing$, then $sh(x) = sh[xv_3v_7] + sh[v_1v_5] + sh[v_2v_4v_6v_8] \le \sigma_{533} +\sigma_{55} + \sigma_{3334} \le \frac{38}{15}$, and we are done. 
Thus, we assume $\{v_{11}, v_{15}, v_{19}, v_{13}\} \cap S = \varnothing$, which leads to the second problem configuration, which will be handled separately.

\begin{figure}[ht]
    \centering
    \begin{tabular}{c@{\hspace{4em}}c}
    \includegraphics[width=0.2\textwidth]{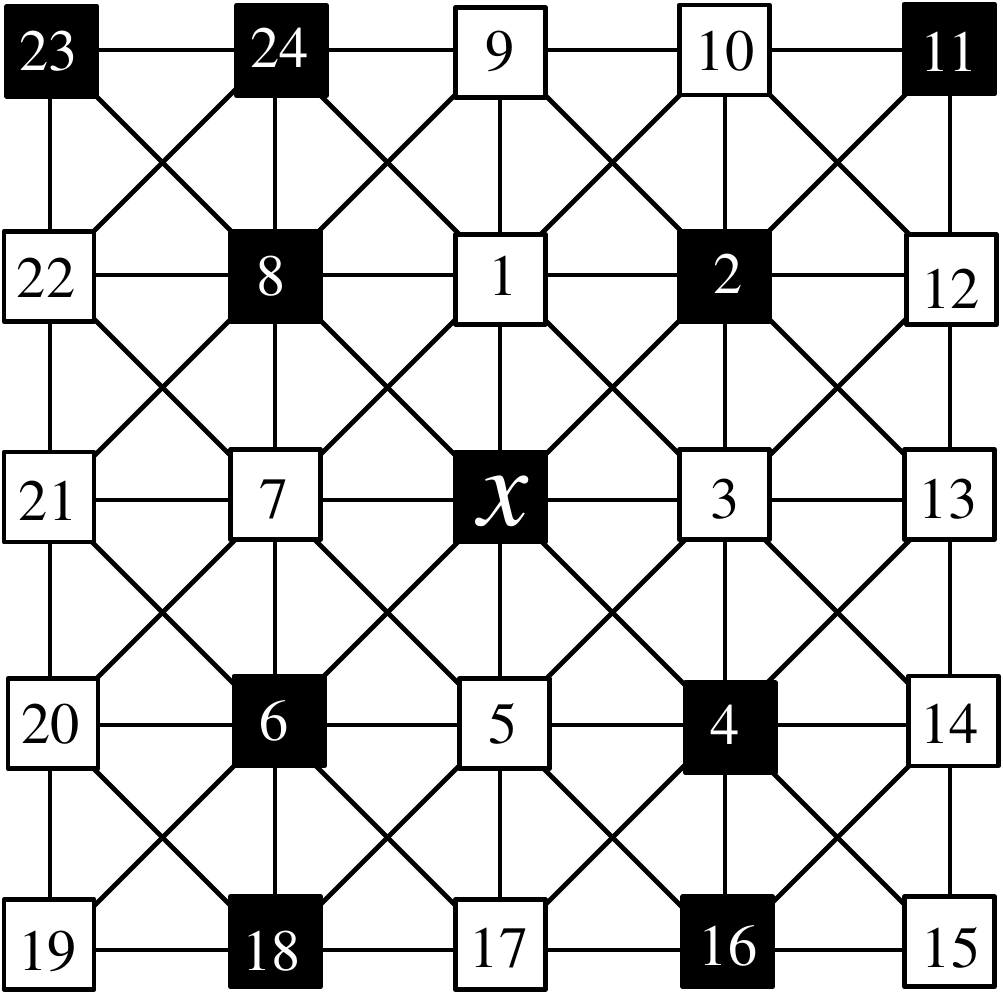} & \includegraphics[width=0.2\textwidth]{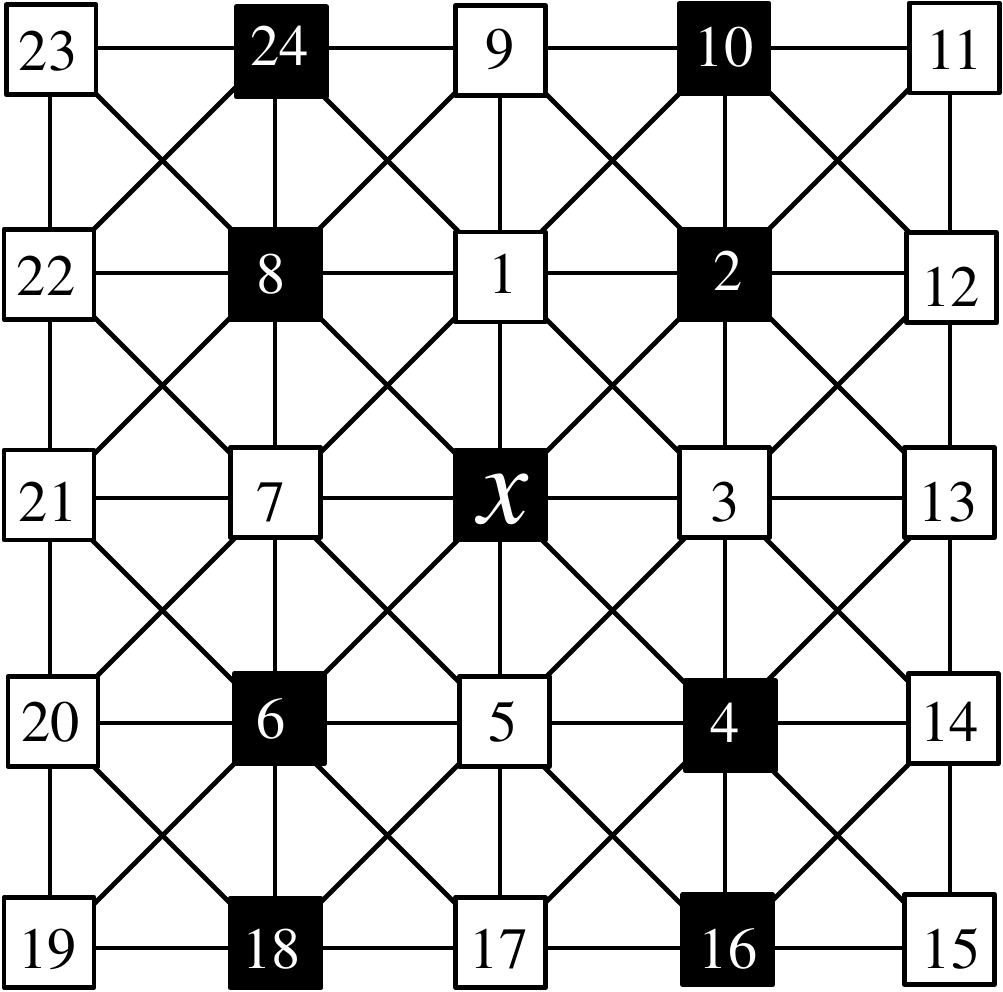} \\
    Problem Case 1 & Problem Case 2
    \end{tabular}
    \caption{X-shape problem cases}
    \label{fig:xshape-probs}
\end{figure}

In either of the Problem Cases---as shown in Figure~\ref{fig:xshape-probs}---we see that $sh(x) \le \max \left\{\frac{77}{30},\frac{13}{5}\right\} = \frac{13}{5}$.
To handle these problems, we consider the maximum share of $v_4$, which will be the same in either configuration.
We observe that we need $dom(v_{17}) \ge 6$ to distinguish from $v_5$, $dom(v_{13}) \ge 4$ to distinguish from $v_3$, and $sh[v_{14}v_{15}] \le \sigma_{34}$ to distinguish $v_{14}$ from $v_{15}$.
Thus $sh(x) = sh[xv_3v_4v_5v_{17}] + sh[v_{13}v_{14}v_{15}v_{16}] = \sigma_{53356} + \sigma_{4343} = \frac{12}{5} < \frac{38}{15}$.
Additionally, the only detector vertex adjacent to $v_4$ other than $x$ is $v_{16}$, which previous cases have shown to have share at most $\frac{38}{15}$.
Thus, $v_4$ can safely accept $\frac{38}{15} - \frac{12}{5} = \frac{2}{15}$ from $x$.
Therefore, $\widehat{sh}(x) \le \frac{13}{5} - \frac{2}{15} = \frac{37}{15} < \frac{38}{15}$ and we are done.

\begin{wrapfigure}{r}{0.2\textwidth}
    \centering
    \includegraphics[width=0.2\textwidth]{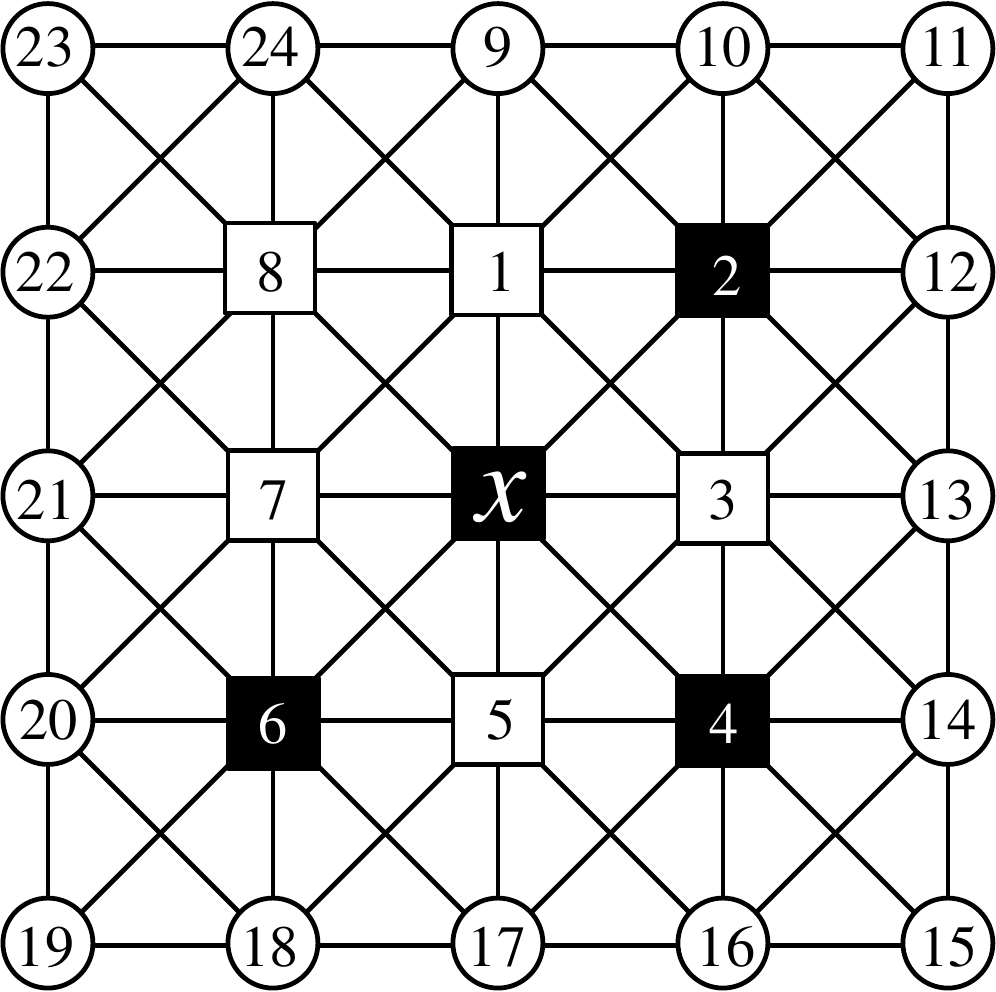}
    \caption{$\lambda$}
    \label{fig:err-ld-case-lambda}
\end{wrapfigure}
\newcase{Lambda}
Consider the case where $x$ is dominated by three neighbors in a ``lambda'' shape, as in Figure \ref{fig:err-ld-case-lambda}.
We will assume that $\{v_1,v_3,v_5,v_7,v_8\} \cap S = \varnothing$, as otherwise we fall into previous cases.

\begin{wrapfigure}[8]{r}{0.25\textwidth}
    \centering
    \vspace{-2em}
    \includegraphics[width=0.2\textwidth]{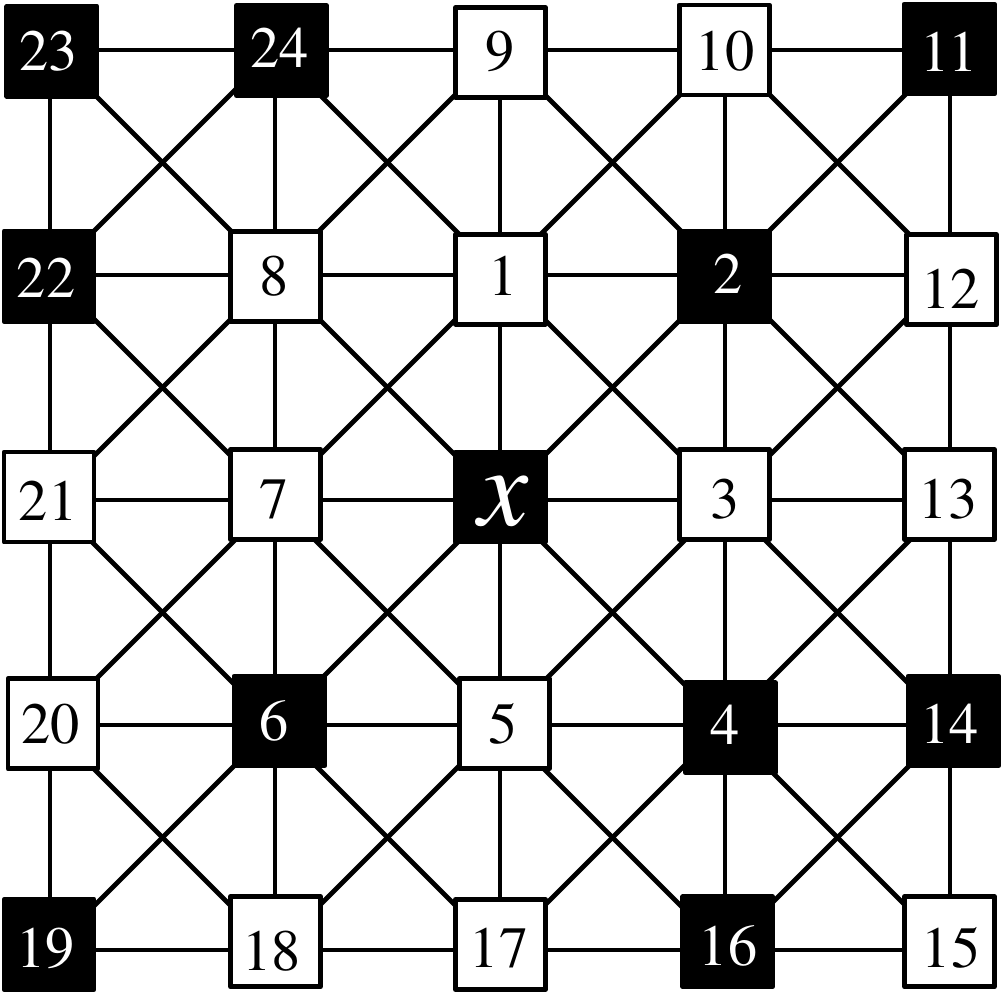}
    \caption{$\lambda$ problem case}
    \label{fig:lambda-prob}
\end{wrapfigure}

Firstly, we observe that $\{v_3,v_5\} \subseteq D_{4+}$ to be distinguish from $x$; so $sh[xv_3v_5] \le \sigma_{444}$.
Suppose $\{v_{20},v_{21}\} \cap S \neq \varnothing$, then we find $sh[v_6v_7] \le \sigma_{35}$ to distinguish $\{v_6,v_7\}$.
By symmetry, $\{v_{9},v_{10}\} \cap S \neq \varnothing$ implies 
$sh[v_1v_2] \le \sigma_{35}$.
Therefore, if $\{v_{20},v_{21}\} \cap S \neq \varnothing$ and $\{v_9,v_{10}\} \cap S \neq \varnothing$, then $sh(x) \le \sigma_{444} + \sigma_{35} + \sigma_{35} +  \sigma_{33} < \frac{38}{15}$ and we are done; otherwise by symmetry we assume that $\{v_9,v_{10}\} \cap S = \varnothing$, requiring $v_{24} \in S$ to 3-dominate $v_1$.
We then see that $dom(v_8) \ge 4$ to be distinguished from $v_1$.
Suppose $\{v_{20},v_{21}\} \cap S \neq \varnothing$ then $sh[v_6v_7] \le \sigma_{35}$ and we are done; thus, we assume $\{v_{20},v_{21}\} \cap S = \varnothing$, requiring $\{v_{22},v_{23}\} \subseteq S$ to 4-dominate $v_8$.
If $\{v_{12},v_{13}\} \cap S \neq \varnothing$, then $sh[v_2v_3] \le \sigma_{35}$ and we are done; thus, by symmetry, we assume that $\{v_{12},v_{13},v_{17},v_{18}\} \cap S = \varnothing$, requiring $v_{11} \in S$ and $v_{19} \in S$ to 3-dominate $v_2$ and $v_6$, respectively.
Additionally, we require $v_{14} \in S$ and $v_{16} \in S$ to 4-dominate $v_3$ and $v_5$, respectively.
If $v_{15} \in S$ then $dom(v_4) \ge 5$ and we are done, so we assume $v_{15} \notin S$, leading to the configuration shown in Figure~\ref{fig:lambda-prob}.

To handle this problem case, we look at the share of $v_4$. 
We observe that we need $dom(v_{15}) \ge 4$ and $dom(v_{13}) \ge 5$ to be distinguished from $v_4$ and $v_3$, respectively.
Therefore, $sh(v_4) \le \frac{2}{5} +\frac{5}{4} +\frac{2}{3} = \frac{139}{60}$. 
Additionally, $v_4$ is only adjacent to  $v_{14}$ and $v_{16}$, which preceding Cases have proven to have share at most $\frac{38}{15}$.
Thus, $v_4$ has no neighboring detectors with share exceeding $\frac{38}{15}$ other than $x$ itself, so it can safely accept $\frac{38}{15} - \frac{139}{60} = \frac{13}{60}$ discharge from $x$.
Therefore, $\widehat{sh}(x) \le \frac{31}{12} - \frac{13}{60} = \frac{71}{30} < \frac{38}{15}$ and we are done.

\begin{wrapfigure}{r}{0.2\textwidth}
    \centering
    \vspace{-1em}
    \includegraphics[width=0.2\textwidth]{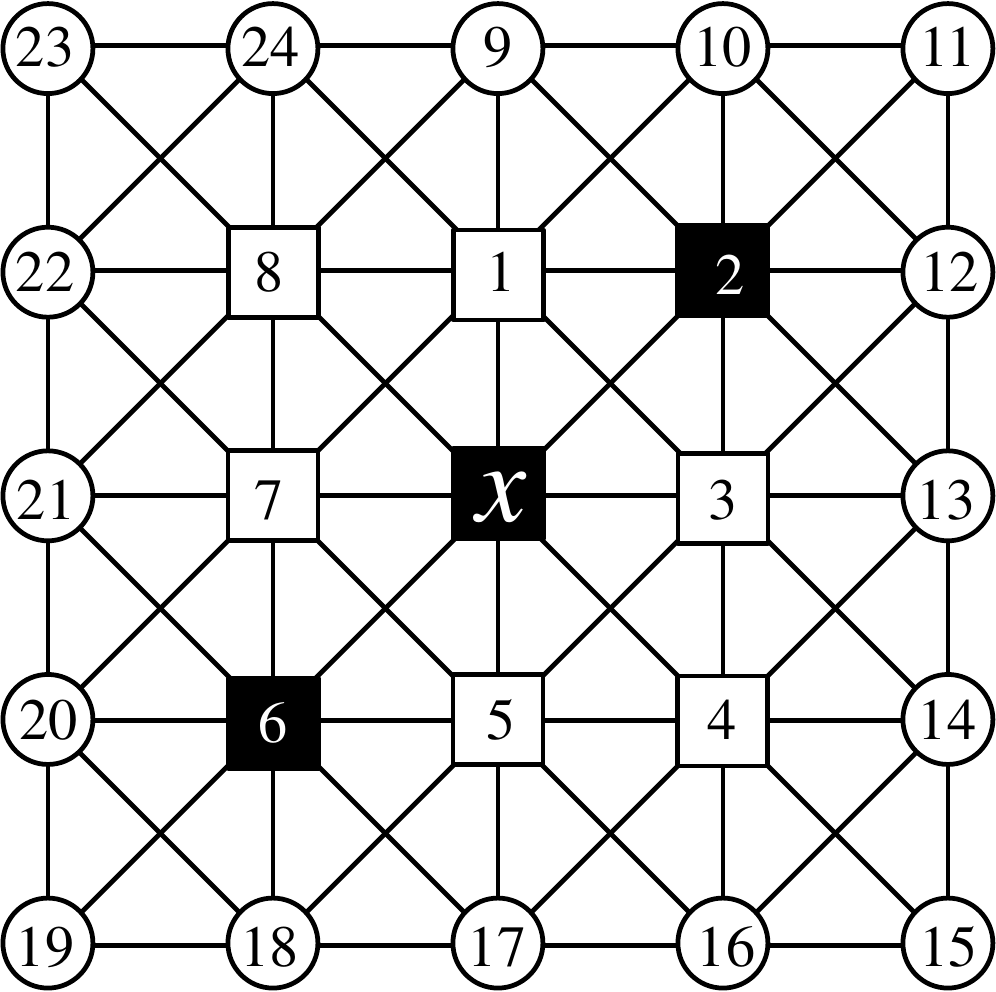}
    \caption{Diagonal}
    \label{fig:err-ld-case-diagonal}
\end{wrapfigure}

\vspace{1em}
\newcase{Diagonal}
Consider the case where $x$ is dominated by two neighbors in a ``diagonal'' shape, as in Figure \ref{fig:err-ld-case-diagonal}.
We can assume that $\{v_1, v_3, v_4,v_5,v_7, v_8\} \cap S = \varnothing$, as otherwise we fall into previous cases.
As this case is perhaps the most complicated, we will begin by proving two smaller claims:

\begin{claim}\label{claim:34-34}
$|\{w \in A : dom(w) = 3\}| \le 1$ for $A = \{v_1,v_3\}$ or $A = \{v_5,v_7\}$.
\end{claim}
\begin{proof}
By symmetry, we need only examine one event: suppose that $dom(v_1) = dom(v_3) = 3$; then $v_1$ and $v_3$ are not distinguished, a contradiction.
\end{proof}

\begin{claim}\label{claim:1k2k}
If $dom(v_1) \ge k$ or $dom(v_3) \ge k$ in the diagonal configuration, then $dom(v_2) \ge k$.
\end{claim}
\begin{proof}
By Theorem~\ref{theo:err-ld-char}, we know $k \le 3$ is always satisfied, and vertices $v_1$ and $v_3$ can be at most 5-dominated, so we consider two cases: $k = 4$ and $k = 5$.

For $k = 4$, without loss of generality assume to the contrary that $dom(v_3) \ge 4$ but $dom(v_2) \le 3$.
We see that $v_2$ is already dominated by itself and $x$, and we consider the third dominator for $v_2$.
If $\{v_9, v_{10}, v_{11}\} \cap S \neq \varnothing$, then $dom(v_3) < 4$, a contradiction; thus, we assume $\{v_9, v_{10}, v_{11}\} \cap S = \varnothing$.
We see that $v_2$ and $v_3$ are not distinguished, a contradiction.

For $k = 5$, without loss of generality, assume to the contrary that $dom(v_3) \ge 5$ but $dom(v_2) \le 4$.
We require $\{v_{12}, v_{13}, v_{14}\} \subseteq S$ to 5-dominate $v_3$; thus $dom(v_2) \ge 4$, so $dom(v_2) = 4$ and we observe that $v_2$ and $v_3$ are not distinguished, a contradiction.
\end{proof}

By symmetry, Claim~\ref{claim:1k2k} also applies to $dom(v_5)$ and $dom(v_7)$ influencing $dom(v_6)$.
\vspace{1em}

Resuming with the diagonal case expansion, suppose $\{v_9,v_{10}\} \cap S = \varnothing$, then $dom(v_1) = 3$ with $v_{24} \in S$.
Then Claim~\ref{claim:34-34} yields that $dom(v_3) \ge 4$, and Claim~\ref{claim:1k2k} yields that $dom(v_2) \ge 4$.
Additionally, $\{v_5,v_7\} \cap D_{4+} \neq \varnothing$, implying $dom(v_6) \ge 4$.
We also require $dom(v_8) \ge 4$ to be distinguished from $v_1$.
We now have five 4-dominated vertices; if we have any more or if any of them is at least 5-dominated then we are done.
Thus, we will assume $\{v_1,v_4\} \subseteq D_3$, $\{v_2,v_3,v_6,v_8\} \subseteq D_4$, and $(dom(v_5),dom(v_7)) \in \{(3,4),(4,3)\}$.
If $\{v_{17},v_{18}\} \subseteq S$ then $v_5$ and $v_6$ are not distinguished, a contradiction.
If $\{v_{17},v_{18}\} \cap S = \varnothing$, then $v_4$ and $v_5$ are not distinguished, a contradiction.
Therefore, we assume $|\{v_{17},v_{18}\} \cap S| = 1$; then $dom(v_5) = 4$ to be distinguished from $v_6$, requiring $dom(v_7) = 3$.
If $v_{20} \in S$ then $v_6$ and $v_7$ are not distinguished, a contradiction, so we assume $v_{20} \notin S$; similarly, we can assume $v_{21} \notin S$, requiring $\{v_{22},v_{23}\} \subseteq S$.
If $\{v_{12},v_{13}\} \subseteq S$ then $v_2$ and $v_3$ are not distinguished, a contradiction.
If $\{v_{12},v_{13}\} \cap S = \varnothing$ then $dom(v_2) = 3$, a contradiction.
Therefore, we assume $|\{v_{12},v_{13}\} \cap S| = 1$, providing symmetry with $\{v_{17},v_{18}\}$.
If $v_{18} \notin S$ then $v_4$ and $v_5$ are not distinguished, a contradiction, so we assume $v_{18} \in S$, requiring $v_{17} \notin S$; by symmetry $v_{12} \in S$ and $v_{13} \notin S$.
Then we require $\{v_{16},v_{19}\} \subseteq S$ and $\{v_{11},v_{14}\} \subseteq S$.
Finally, $v_{15} \notin S$ because $dom(v_4) = 3$.
This gives us Problem Case 1, which will be handled separately.

Otherwise, $\{v_9,v_{10}\} \cap S \neq \varnothing$, and by symmetry $\forall A \in \{\{v_{12},v_{13}\}, \{v_{17},v_{18}\}, \{v_{20},v_{21}\}\}$, $A \cap S \neq \varnothing$.
Consider the corner $\{v_1,v_2,v_3\}$: we now know that $sh[v_1v_2] \le \max\{\sigma_{35},\sigma_{44}\}$, and similarly for $sh[v_2v_3]$.
Thus, $sh[v_1v_2v_3] \le \max\{\sigma_{345},\sigma_{444}\} = \sigma_{345}$.
Suppose $sh[v_1v_2v_3] \le \sigma_{444}$, then by symmetry $sh[v_5v_6v_7] \le \sigma_{345}$ and we are done.
So we can assume $\exists w \in \{v_1,v_2,v_3\}$ such that $dom(w) = 3$, and by Claims \ref{claim:34-34} and \ref{claim:1k2k}, $w$ is unique.
Suppose $w = v_2$, then by Claim~\ref{claim:1k2k} we reach a contradiction; without loss of generality let $w = v_1$; then $dom(v_2) \ge 5$ to distinguish from $v_1$ because $\{v_9,v_{10}\} \cap S \neq \varnothing$.
If $dom(v_3) \ge 5$ then we are done, so we assume $dom(v_3) = 4$.
If $dom(v_2) \ge 6$ we are done, so we assume $dom(v_2) = 5$.
With symmetry we can also assume that $dom(v_6) = 5$, that $(dom(v_5),dom(v_7)) \in \{(3,4), (4,3)\}$, and that $dom(v_4) = dom(v_8) = 3$.
If $v_9 \in S$ or $v_{24} \in S$ then $v_1$ and $v_8$ are not distinguished, a contradiction; thus, $v_9 \notin S$ and $v_{24} \notin S$, requiring $v_{10} \in S$.
If $dom(v_7) = 3$ then by symmetry $v_{21} \notin S$ and $v_{22} \notin S$, but we see that $v_8$ cannot be 3-dominated, a contradiction; therefore $dom(v_7) = 4$, meaning $dom(v_5) = 3$.
By symmetry with $\{v_9,v_{10},v_{24}\}$, we can assume $\{v_{16},v_{17},v_{18}\} \cap S = \{v_{18}\}$.
If $v_{21} \in S$ then $v_{20} \in S$ to distinguish $v_7$ from $v_8$, however this results in $v_6$ and $v_7$ not being distinguished, a contradiction; thus, we assume $v_{21} \notin S$.
Then to 5-dominate $v_6$ and 3-dominate $v_8$ we require $\{v_{19},v_{20},v_{22},v_{23}\} \subseteq S$.
By symmetry, we can also assume that $v_{13} \notin S$ and $\{v_{11},v_{12},v_{14},v_{15}\} \subseteq S$.
This results in Problem Case 2, which will be handled separately.

\begin{figure}[ht]
    \centering
    \begin{tabular}{c@{\hspace{5em}}c}
        \includegraphics[width=0.2\textwidth]{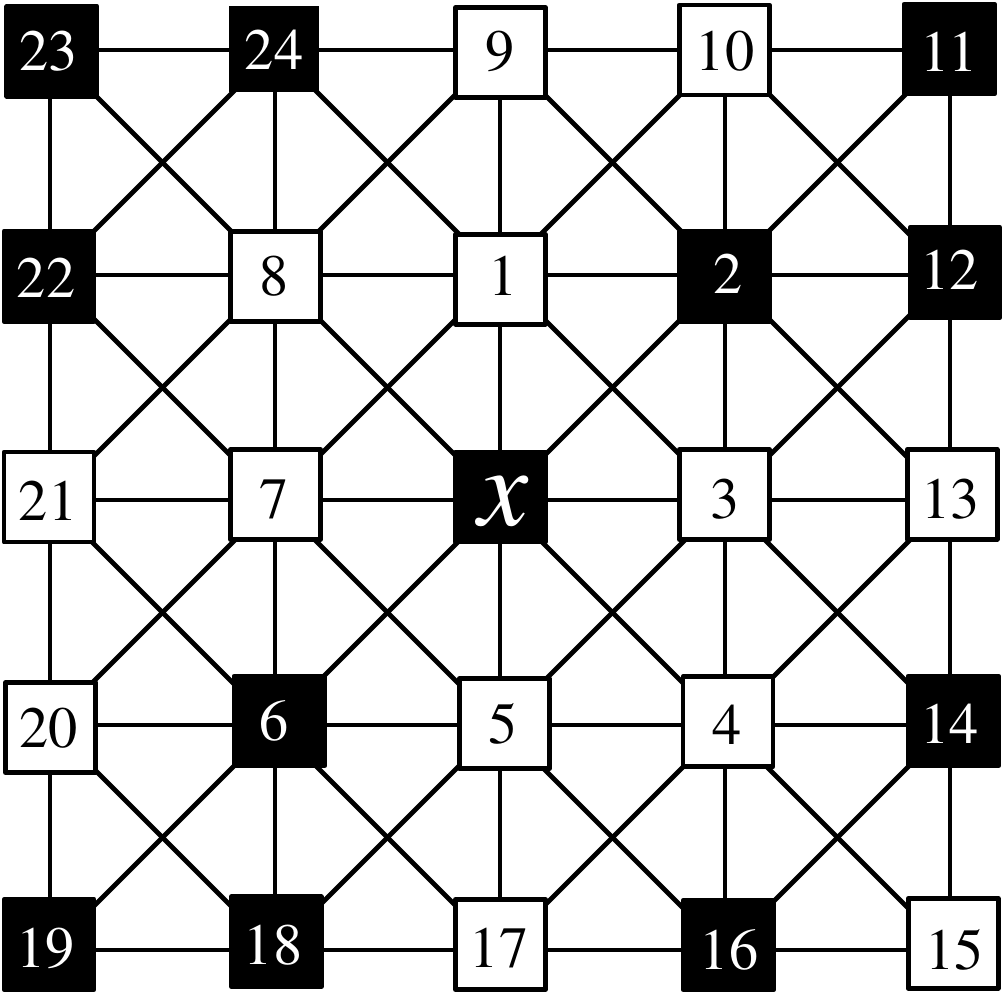} & \includegraphics[width=0.2\textwidth]{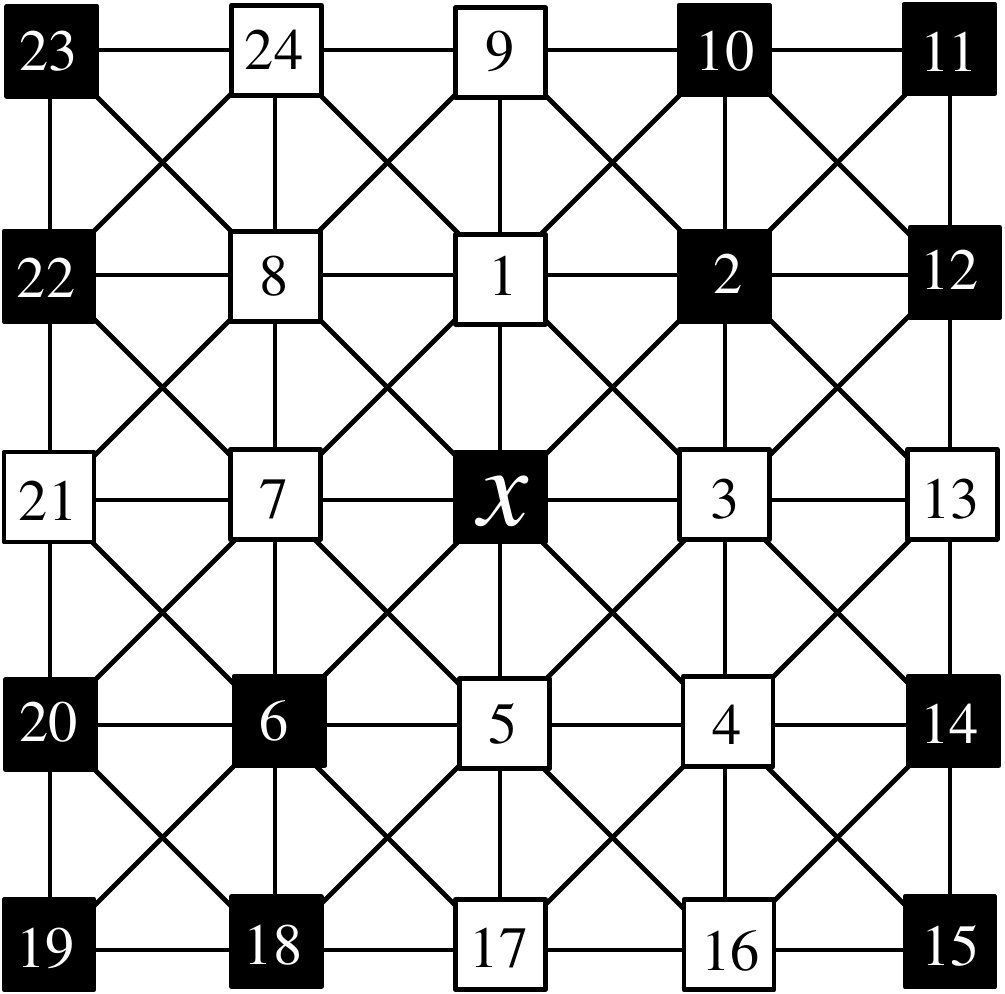} \\ Problem Case 1 & Problem Case 2
    \end{tabular}
    \caption{Diagonal problem cases}
    \label{fig:err-ld-diagonal-prob}
\end{figure}

To handle the first ``diagonal'' problem case, we consider the share of $v_6$.  
We observe that we need $dom(v_{17}) \ge 5$, $dom(v_{20}) \ge 4$, and $dom(v_{21}) \ge 4$ to distinguish from $v_5$, $v_6$, and $v_7$, respectively.  
We also see that $\{v_{18},v_{19}\} \cap D_{4+} \neq \varnothing$ to distinguish from each other.
Therefore, $sh(v_6) \le \frac{1}{5} + \frac{5}{4} + \frac{3}{3} = \frac{49}{20}$.
Additionally, the only detectors adjacent to $v_6$ other than $x$ are $v_{18}$ and $v_{19}$, which preceding Cases have proven to have share at most $\frac{38}{15}$.
Thus, $v_6$ has no neighboring detectors with share exceeding $\frac{38}{15}$ other than $x$ itself, so it can safely accept $\frac{38}{15} - \frac{49}{20} = \frac{1}{12}$ discharge from $x$.
Therefore, $\widehat{sh}(x) \le \frac{31}{12} - \frac{1}{12} = \frac{5}{2} < \frac{38}{15}$ and we are done.

To handle the other ``diagonal'' problem case, we consider the share of $v_6$.  
We observe that we need $dom(v_{17}) \ge 4$ and $dom(v_{21}) \ge 5$ to distinguish from $v_7$ and $v_5$, respectively.
We also see that $|\{v_{18},v_{19},v_{20}\} \cap D_{5+}| \ge 2$ to be distinguished. 
Therefore, $sh(v_6) \le \frac{4}{5} +\frac{3}{4} +\frac{2}{3} = \frac{133}{60}$.
Additionally, the only detectors adjacent to $v_6$ other than $x$ are $v_{18}$, $v_{19}$ and $v_{20}$, which preceding Cases have proven to have share at most $\frac{38}{15}$.
Thus, $v_6$ has no neighboring detectors with share exceeding $\frac{38}{15}$ other than $x$ itself, so it can accept $\frac{38}{15} - \frac{133}{60} = \frac{19}{60}$ discharge from $x$.
Therefore, $\widehat{sh}(x) \le \frac{77}{30} - \frac{19}{60} = \frac{9}{4} < \frac{38}{15}$, completing the diagonal configuration.

\begin{wrapfigure}{r}{0.2\textwidth}
    \centering
    \vspace{-1em}
    \includegraphics[width=0.2\textwidth]{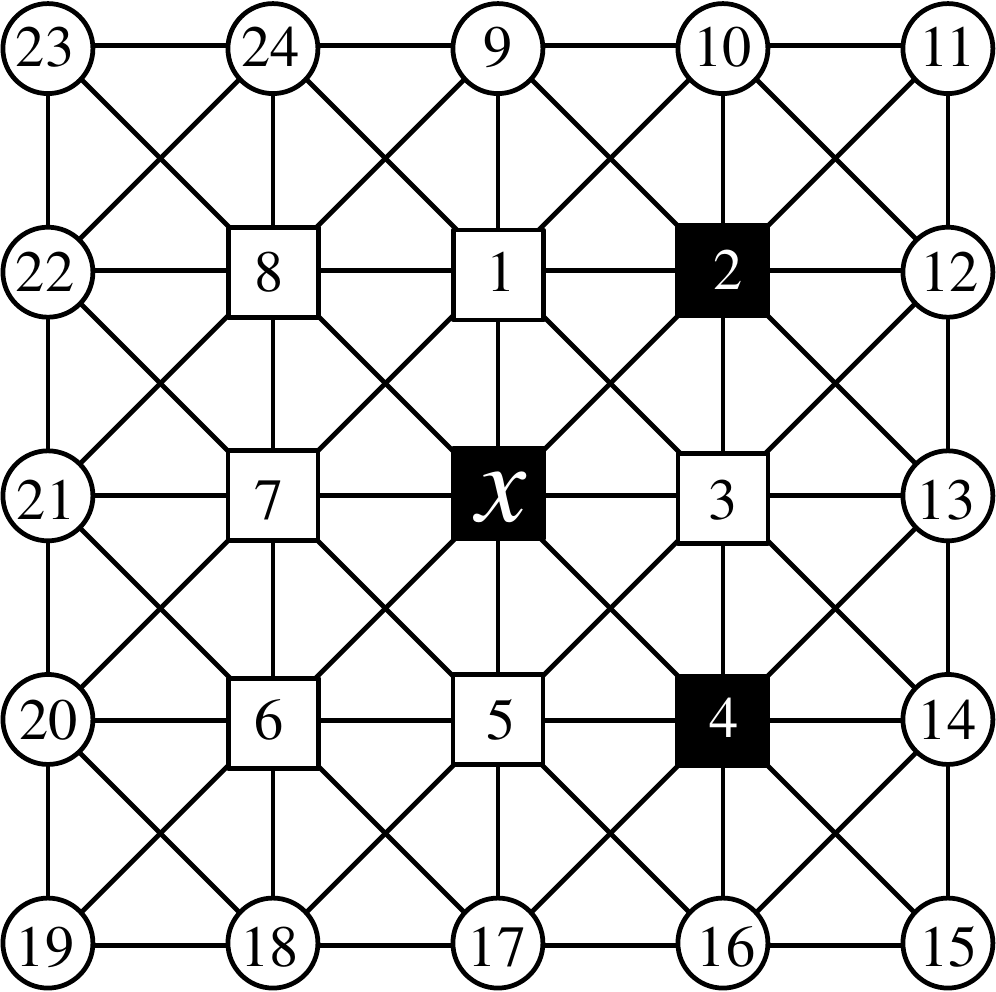}
    \caption{Wedge}
    \label{fig:err-ld-case-wedge}
\end{wrapfigure}

\newcase{Wedge}
As the last case, we consider when $x$ is dominated by two neighbors in a ``wedge'' shape, as in Figure \ref{fig:err-ld-case-wedge}.
We can assume that $\{v_1, v_3, v_5, v_6, v_7, v_8\}  \cap S = \varnothing$, as otherwise we fall into previous cases.

First, we observe that $dom(v_3) \ge 5$ to be distinguished from $x$, and $dom(v_7) \in \{3,4\}$.
If $dom(v_7) = 4$, then $dom(v_6) \ge 5$ is required to be distinguished from $v_7$; by symmetry, $dom(v_8) \ge 5$, so $sh(x) \le \frac{3}{5} + \frac{1}{4} + \frac{5}{3} \le \frac{38}{15}$ and we are done.
If $v_7$ is 3-dominated with $\{v_{21},v_{22}\} \subseteq S$, then we observe that $dom(v_8) \ge 6$ to be distinguished from $v_7$, requiring $\{v_9,v_{23},v_{24}\} \subseteq S$; thus, $dom(v_1) \ge 4$.
We also find that $dom(v_6) \ge 4$ to be distinguished from $v_7$.
Therefore, $sh(x) \le \frac{1}{6} + \frac{1}{5} + \frac{2}{4} + \frac{5}{3} = \frac{38}{15}$ and we are done.
By symmetry, we get the same results if $v_7$ is 3-dominated with $\{v_{20},v_{21}\} \subseteq S$.

Therefore, we assume $v_7$ is 3-dominated with $\{v_{20},v_{22}\} \subseteq S$.
We observe that $dom(v_8) \ge 4$ to be distinguished from $v_7$; by symmetry, this is also true for $v_6$.
Suppose that $v_{19} \notin S$.
Then $\{v_{17},v_{18}\} \subseteq S$ to 4-dominate $v_6$, and the only way to distinguish $v_5$ and $v_6$ is to have $v_{16} \in S$, so $dom(v_5) = 5$ and $dom(v_4) \ge 4$.
Therefore, $sh(x) \le \frac{2}{5} + \frac{3}{4} + \frac{4}{3} = \frac{149}{60} < \frac{38}{15}$ and we are done; otherwise, we assume $v_{19} \in S$, and by symmetry we also assume $v_{23} \in S$.
Next, suppose that $\{v_{18},v_{24}\} \cap S = \varnothing$, then $v_{17} \in S$ and $v_9 \in S$ to 4-dominate $v_6$ and $v_8$, respectively.
We observe that we need $dom(v_2) \ge 5$ to be distinguished from $v_1$, and by symmetry $dom(v_4) \ge 5$ as well; thus, $sh(x) \le \frac{3}{5} + \frac{2}{4} + \frac{4}{3} = \frac{73}{30} < \frac{38}{15}$ and we are done.
Otherwise we can assume that $\{v_{18},v_{24}\} \cap S \neq \varnothing$, and without loss of generality let $v_{18} \in S$; we now consider the two cases depending on whether $v_{24}$ is a detector.

\begin{figure}[b!]
    \centering
    \begin{tabular}{c@{\hspace{5em}}c}
        \includegraphics[width=0.2\textwidth]{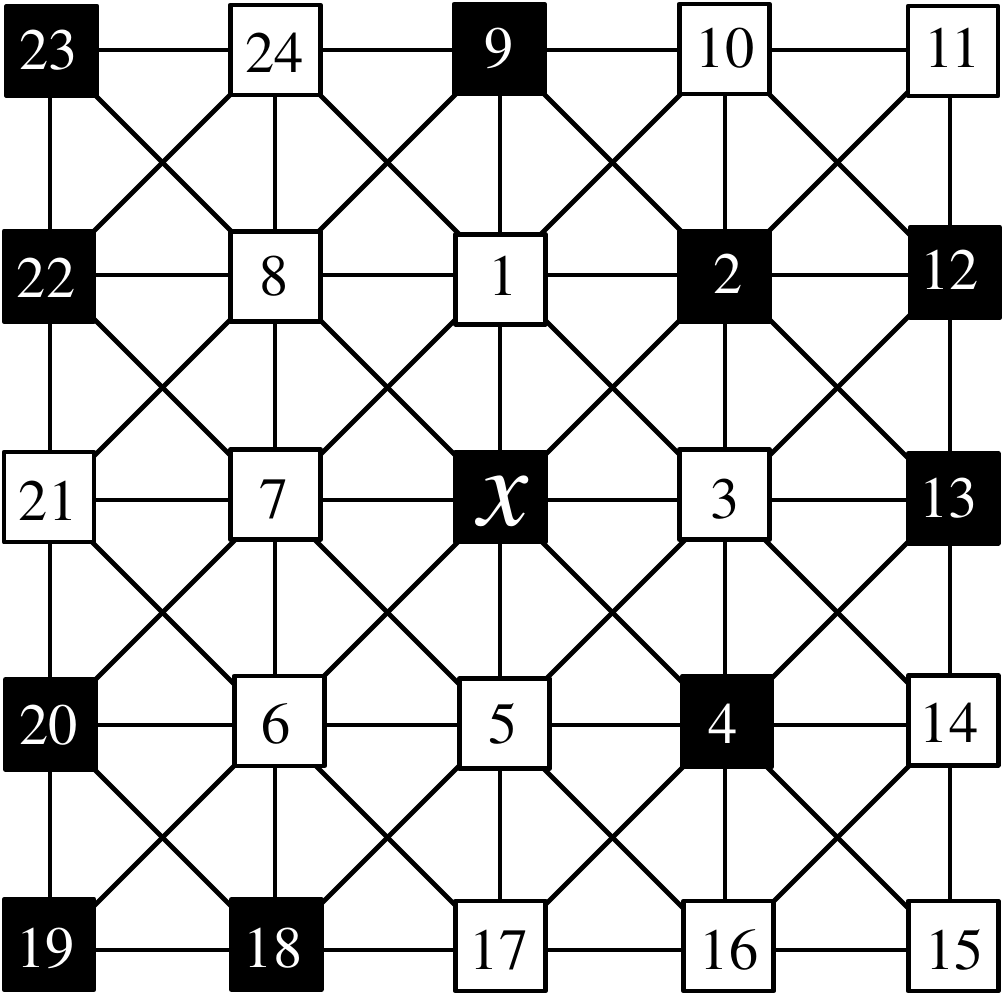} & \includegraphics[width=0.2\textwidth]{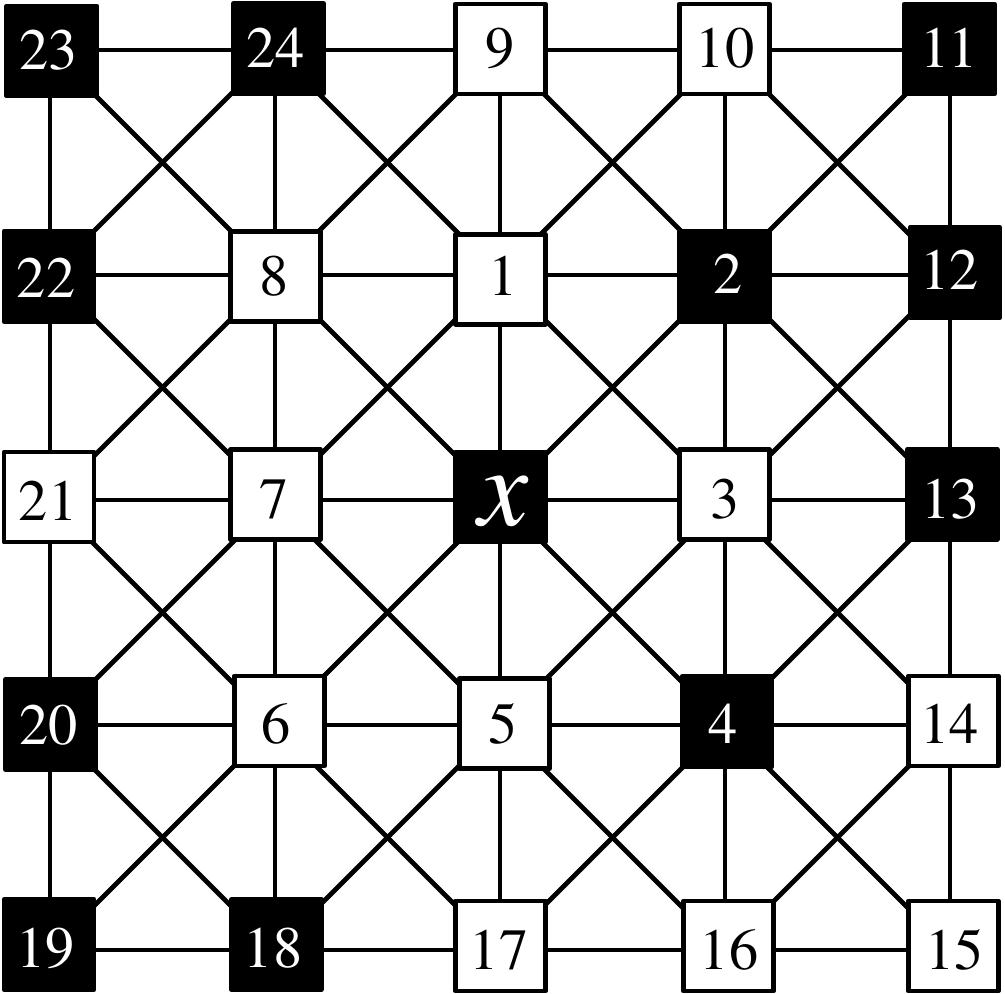} \\
        Problem 1 & Problem 2
    \end{tabular}
    \begin{tabular}{c@{\hspace{5em}}c@{\hspace{5em}}c} \\
        \includegraphics[width=0.2\textwidth]{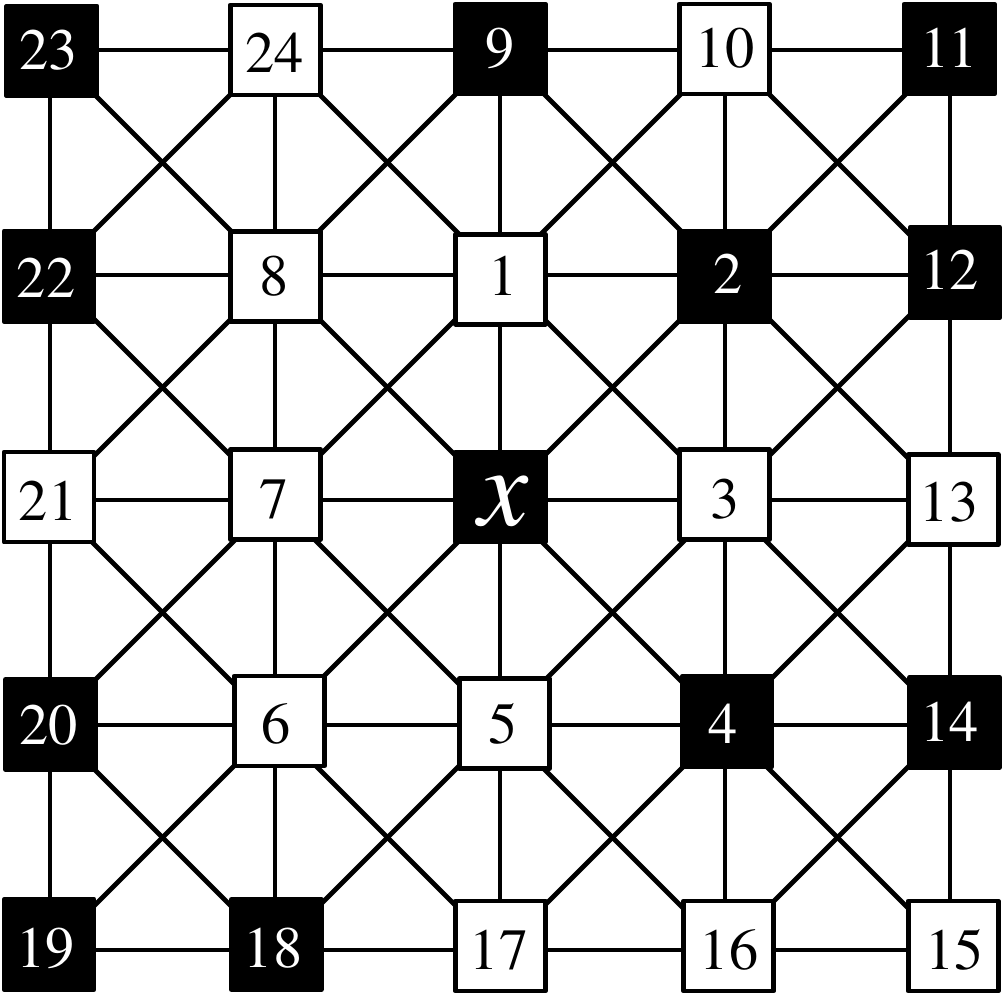} & \includegraphics[width=0.2\textwidth]{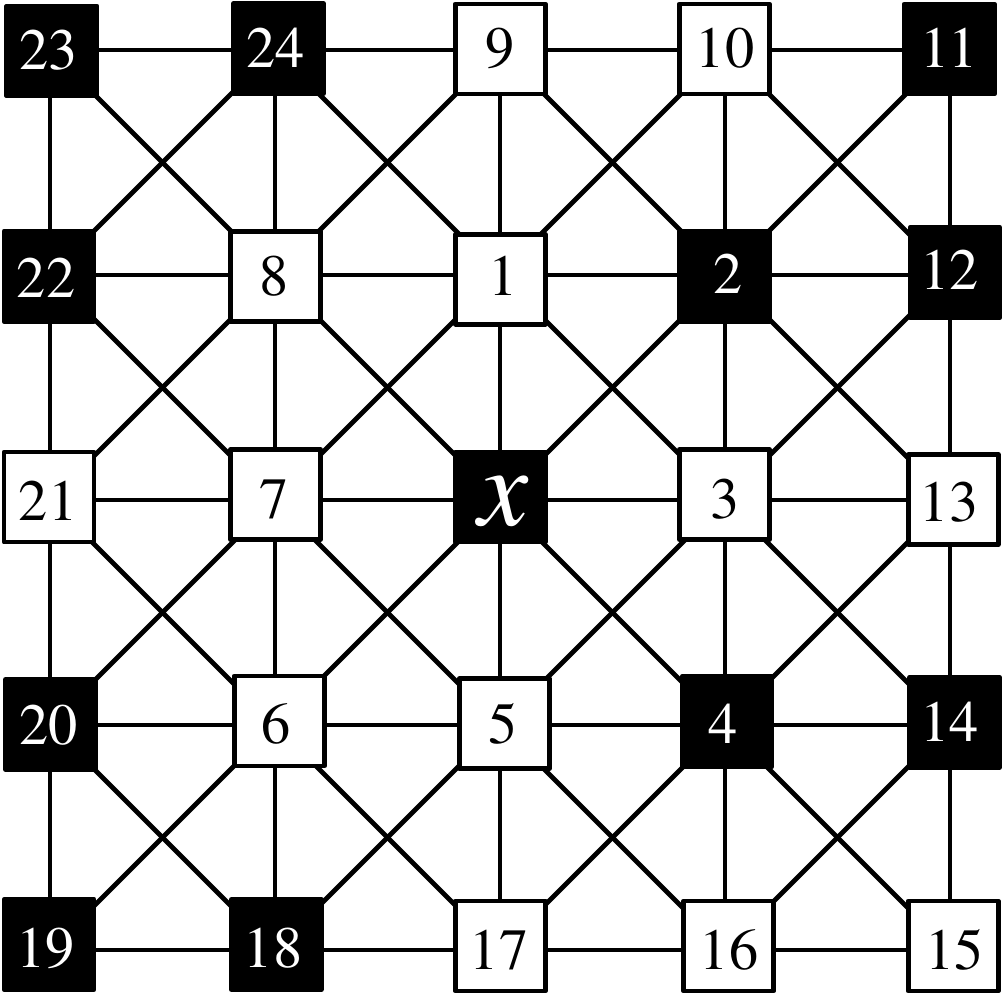} &
        \includegraphics[width=0.2\textwidth]{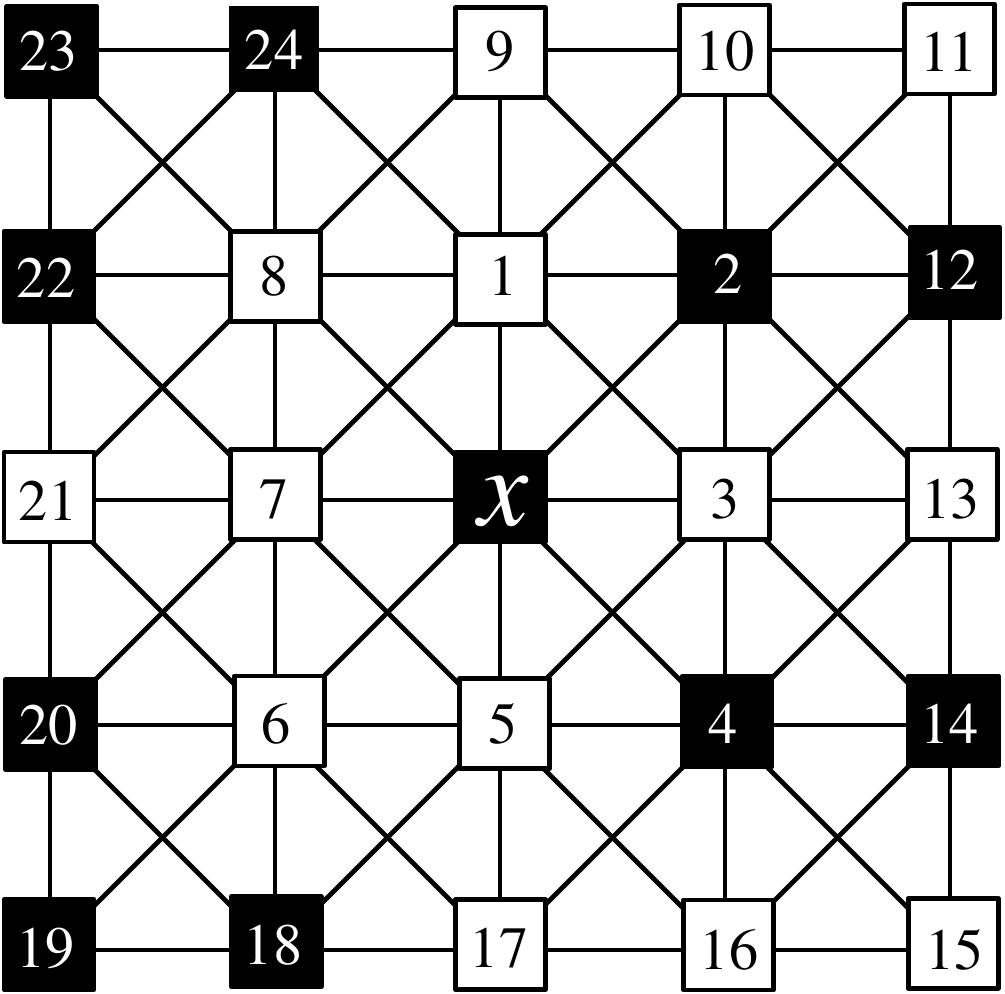} \\
        Problem 3 & Problem 4 & Problem 5
    \end{tabular}
    \caption{Wedge problem cases}
    \label{fig:wedge-probs-345}
\end{figure}

If $v_{24} \notin S$, then $v_9 \in S$ to 4-dominate $v_8$ and we observe that $dom(v_2) \ge 5$ to be distinguished from $v_1$.
If $v_{17} \in S$ then $dom(v_6) = 5$ and $dom(v_5) \ge 4$, resulting in $sh(x) \le \frac{3}{5} + \frac{2}{4} + \frac{4}{3} = \frac{73}{30} < \frac{38}{15}$ and we are done; otherwise we assume $v_{17} \notin S$.
If $\{v_1,v_4,v_5\} \cap D_{4+} \neq \varnothing$ then $sh(x) \le \frac{2}{5} + \frac{3}{4} + \frac{4}{3} = \frac{149}{60} < \frac{38}{15}$ and we are done; thus, we can assume that $\{v_1,v_4,v_5\} \subseteq D_3$.
Therefore, $v_{10} \notin S$ and $v_{16} \notin S$ because $v_1$ and $v_5$ are 3-dominated, respectively.
Suppose $v_{13} \in S$; then $\{v_{14},v_{15}\} \cap S = \varnothing$ because $v_4$ is 3-dominated, so we require $v_{12} \in S$ to 5-dominate $v_3$.
If $v_{11} \in S$ then $dom(v_2) = 6$, resulting in $sh(x) \le \frac{1}{6} + \frac{1}{5} + \frac{2}{4} + \frac{5}{3} = \frac{38}{15}$ and we are done; otherwise $v_{11} \notin S$ and we arrive at Problem Case 1, which will be handled separately.
Otherwise, we consider $v_{13} \notin S$; then $\{v_{12},v_{14}\} \subseteq S$ to 5-dominate $v_3$.
This requires $v_{15} \notin S$ because $v_4$ is 3-dominated, and $v_{11} \in S$ to 5-dominate $v_2$; this leads to Problem Case 3, which will be handled separately.

We now consider when $v_{24} \in S$.
Suppose that $v_9 \in S$; then $dom(v_8) = 5$ and $dom(v_1) \ge 4$.
If $dom(v_2) \ge 4$ then $sh(x) \le \frac{2}{5} + \frac{3}{4} + \frac{4}{3} = \frac{149}{60} < \frac{38}{15}$ and we are done; otherwise, we assume that $dom(v_2) = 3$.
Therefore, $\{v_{10},v_{11},v_{12},v_{13}\} \cap S = \varnothing$, however this results in $v_1$ and $v_2$ not being distinguished, a contradiction.
Therefore, $v_9 \notin S$, and by symmetry $v_{17} \notin S$ as well.
By applying the same argumentation, we can also prove that $v_{10} \notin S$, and by symmetry $v_{16} \notin S$ as well.
Suppose $v_{13} \in S$.
If $dom(v_4) = 3$ then $\{v_{14},v_{15}\} \cap S = \varnothing$, we require $v_{12} \in S$ to 5-dominate $v_3$ and $v_{11} \in S$ to distinguish $v_2$ and $v_3$; this leads to Problem Case 2 and will be handled separately.
Otherwise, $dom(v_4) \ge 4$, and by symmetry $dom(v_2) \ge 4$ as well, so $sh(x) \le \frac{1}{5} + \frac{4}{4} + \frac{4}{3} = \frac{38}{15}$ and we are done.
Therefore, we can assume $v_{13} \notin S$; then $\{v_{12},v_{14}\} \subseteq S$ to 5-dominate $v_3$.
If $\{v_{11},v_{15}\} \subseteq S$ then $\{v_2,v_4\} \subseteq D_{4+}$, so we are done; otherwise, without loss of generality we assume $v_{15} \notin S$.
If $v_{11} \in S$ then we arrive at Problem Case 4, otherwise we arrive at Problem Case 5; both of these will be handled separately later.

\FloatBarrier
\begin{wrapfigure}{r}{0.225\textwidth}
    \centering
    \vspace{-2em}
    \includegraphics[width=0.225\textwidth]{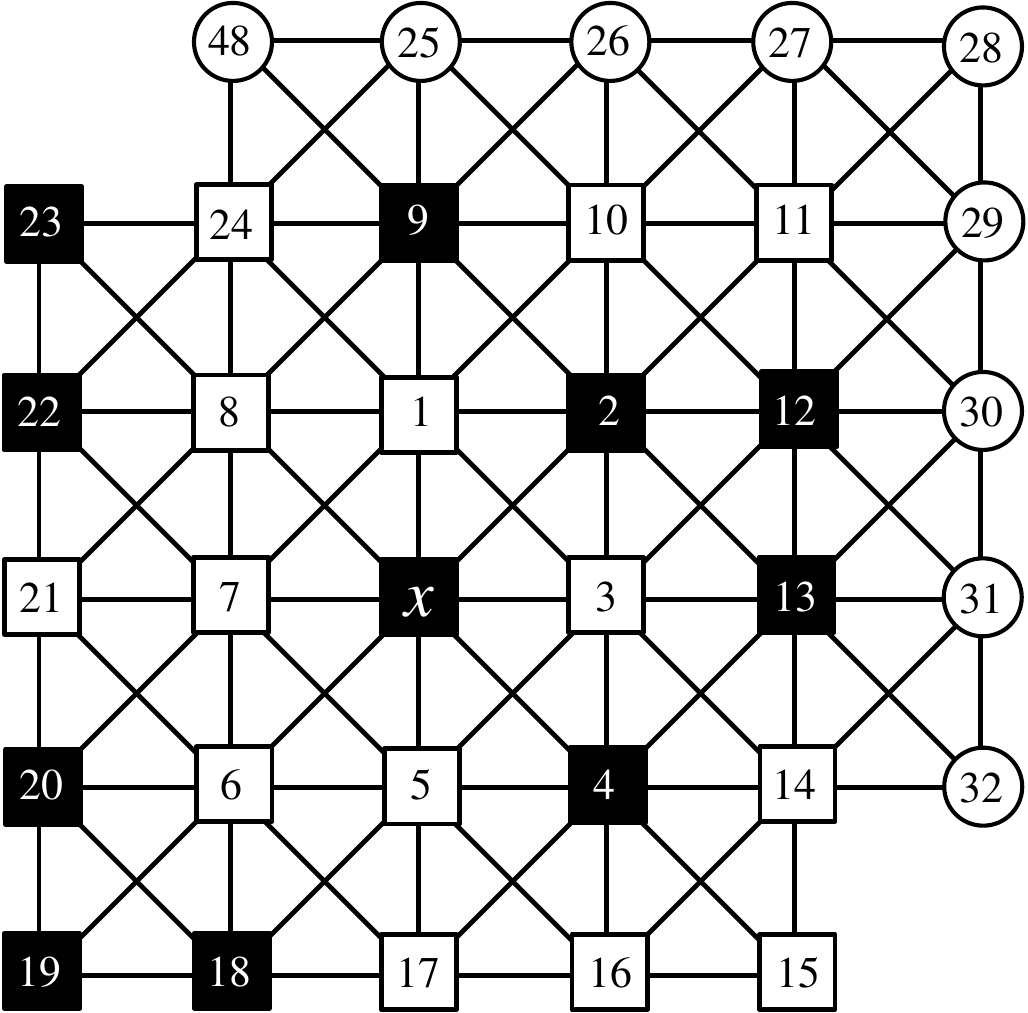}
    \caption{Problem case 1, expanded}
    \label{fig:wedge-prob-1}
\end{wrapfigure}

\textbf{Wedge Problem Case 1:}
To handle Wedge Problem Case 1 configuration, we consider the maximum share of $v_2$; we will show that $sh(v_2) \le \frac{37}{15}$.
We observe that $dom(v_{13}) \ge 5$ and $dom(v_{10}) \ge 4$ to distinguish from $v_3$ and $v_1$, respectively.
If $dom(v_{10}) \ge 5$ then $sh(v_2) \le \frac{4}{5} + \frac{5}{3} = \frac{37}{15}$ and we are done; otherwise we assume $dom(v_{10}) = 4$.
If $\{v_9,v_{11}\} \cap D_{4+} \neq \varnothing$ then we are done; so we assume $\{v_9,v_{11}\} \subseteq D_3$.
If $v_{26} \in S$ or $v_{27} \in S$ to 4-dominate $v_{10}$, then $v_{10}$ and $v_{11}$ are not distinguished, a contradiction.
Thus, $v_{25} \in S$, $v_{26} \notin S$, and $v_{27} \notin S$, but we find that $v_9$ and $v_{10}$ are not distinguished, a contradiction; therefore, $sh(v_2) \le \frac{37}{15}$.
Additionally, $v_2$ is only adjacent to three detectors other than $x$, $v_9$, $v_{12}$, and $v_{13}$.
From previous cases, we know that $sh(v_{12}) \le \frac{38}{15}$ and $sh(v_{13}) \le \frac{38}{15}$; however, $sh(v_9)$ might exceed $\frac{38}{15}$.
Therefore, we conservatively allow $v_2$ to accept discharge from two sources---$x$ and potentially $v_9$---up to $\frac{1}{2} \left[ \frac{38}{15} - \frac{37}{15} \right] = \frac{1}{30}$ from each.
Thus, $\widehat{sh}(x) = \frac{77}{30} - \frac{1}{30} = \frac{38}{15}$ and we are done.

\begin{wrapfigure}[11]{r}{0.225\textwidth}
    \centering
    \vspace{-1em}
    \includegraphics[width=0.225\textwidth]{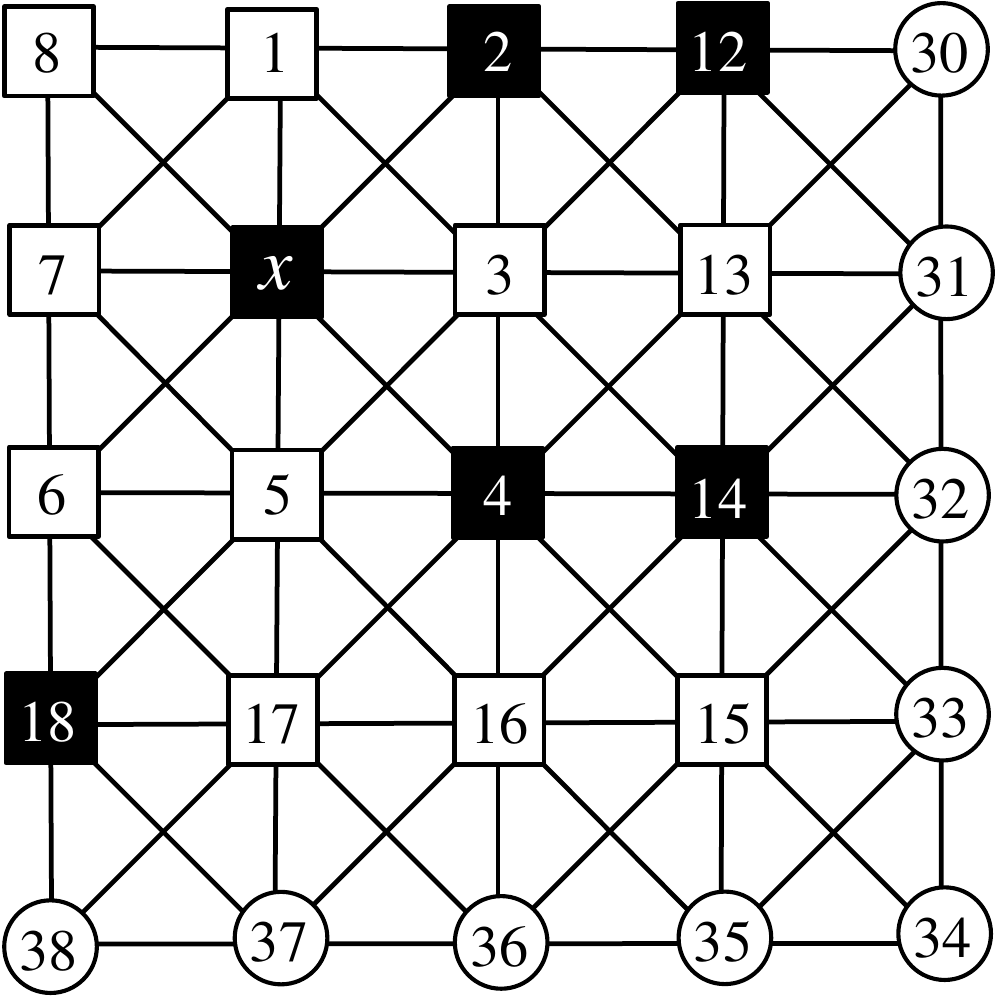}
    \caption{Problems 3--5 common region}
    \label{fig:wedge-prob-5}
\end{wrapfigure}

\textbf{Wedge Problem Case 2:}
To handle this situation, we look at the share of $v_2$. 
We observe that $sh[v_{10}v_{11}] \le \sigma_{35}$ to distinguish $v_{10}$ and $v_{11}$.
Therefore, $sh(v_2) \le \frac{3}{5} +\frac{2}{4}+\frac{4}{3} = \frac{73}{30}$. 
Additionally, $v_2$ is only adjacent to  $v_{11}$, $v_{12}$, $v_{13}$, which ``L'' shape and ``triangle'' Cases have proven to have share at most $\frac{38}{15}$.
Thus, $v_2$ has no neighboring detectors with share exceeding $\frac{38}{15}$ other than $x$ itself, so it can accept $\frac{38}{15} - \frac{73}{30} = \frac{1}{10}$ discharge from $x$.
Therefore, $\widehat{sh}(x) \le \frac{77}{30} - \frac{1}{10} = \frac{37}{15} < \frac{38}{15}$ and we are done.

\textbf{Wedge Problem Cases 3, 4, and 5:}
Consider Problem Cases 3--5, as shown in Figure~\ref{fig:wedge-probs-345}; to handle these cases, we will show that $sh(v_4) \le \frac{49}{20}$.
We observe that $dom(v_{13}) \ge 6$ to distinguish from $v_3$, $dom(v_{17}) \ge 4$ to distinguish from $v_5$.
If $\{v_{35},v_{36}\} \cap S \neq \varnothing$ then we require $sh[v_{15}v_{16}] \le \max \left\{ \sigma_{36},\sigma_{45} \right\} = \sigma_{36}$ to distinguish $v_{15}$ from $v_{16}$, so $sh(v_4) \le \frac{49}{20}$ and we are done; thus, we assume $\{v_{35},v_{36}\} \cap S = \varnothing$, meaning $dom(v_{16}) = 3$, so $dom(v_{15}) \ge 4$.
If $dom(v_{14}) \ge 4$ then we will be done, so we assume $dom(v_{14}) = 3$.
If $v_{30} \notin S$ or $v_{33} \in S$ then we contradict that $dom(v_{13}) \ge 6$, so we assume $v_{30} \in S$ and $v_{33} \notin S$.
Thus, we require $\{v_{32},v_{34}\} \subseteq S$ to 4-dominate $v_{15}$, and $v_{31} \notin S$ because $dom(v_{14}) = 3$.
We see that $v_{14}$ and $v_{15}$ are not distinguished, a contradiction.
Therefore, we have proven that $sh(x_4) \le \frac{49}{20}$.
Additionally, the only detector adjacent to $v_4$ other than $x$ is $v_{14}$, which previous cases have shown to have share no more than $\frac{38}{15}$; so $v_4$ can accept $\frac{38}{15} - \frac{49}{20} = \frac{1}{12}$ discharge from $x$.
Consider Problem Cases 3 and 4: $sh(x) \le \max \left\{ \frac{157}{60}, \frac{77}{30} \right\} = \frac{157}{60}$. Then $\widehat{sh}(x) \le \frac{157}{60} - \frac{1}{12} = \frac{38}{15}$ and we are done.
Consider Problem Case 5: $sh(x) = \frac{27}{10}$.
By symmetry, we can also discharge up to $\frac{1}{12}$ from $x$ into $v_2$, so $\widehat{sh}(x) \le \frac{27}{10} - 2 \times \frac{1}{12} = \frac{38}{15}$, completing the proof of Theorem~\ref{theo:sh-err-ld-king}.
\cendproof


\newpage
\begin{theorem}\label{theo:red-sh}
Let $S \subset V(G)$ be a RED:LD set on $\textrm{K}$. Then, the average share of all vertices in $S$ is at most $\frac{11}{3}$.
\end{theorem}
\cbeginproof
Let $S \subseteq V(G)$ be a RED:LD set on $\textrm{K}$, let the notations $sh$, $dom$, and $k$-dominated implicitly use $S$, and let a vertex pair being \emph{distinguished} denote satisfying the relevant pair property of Theorem~\ref{theo:red-ld-char}.
Let $x \in S$ be a detector vertex; Theorem~\ref{theo:red-ld-char} yields that $|N[x] \cap S| \ge 2$.
There are two non-isomorphic ways to 2-dominate $x$, either ``vertically'' or ``diagonally.''

\begin{wrapfigure}[11]{r}{0.2\textwidth}
    \centering
    \vspace{-1em}
    \includegraphics[width=0.2\textwidth]{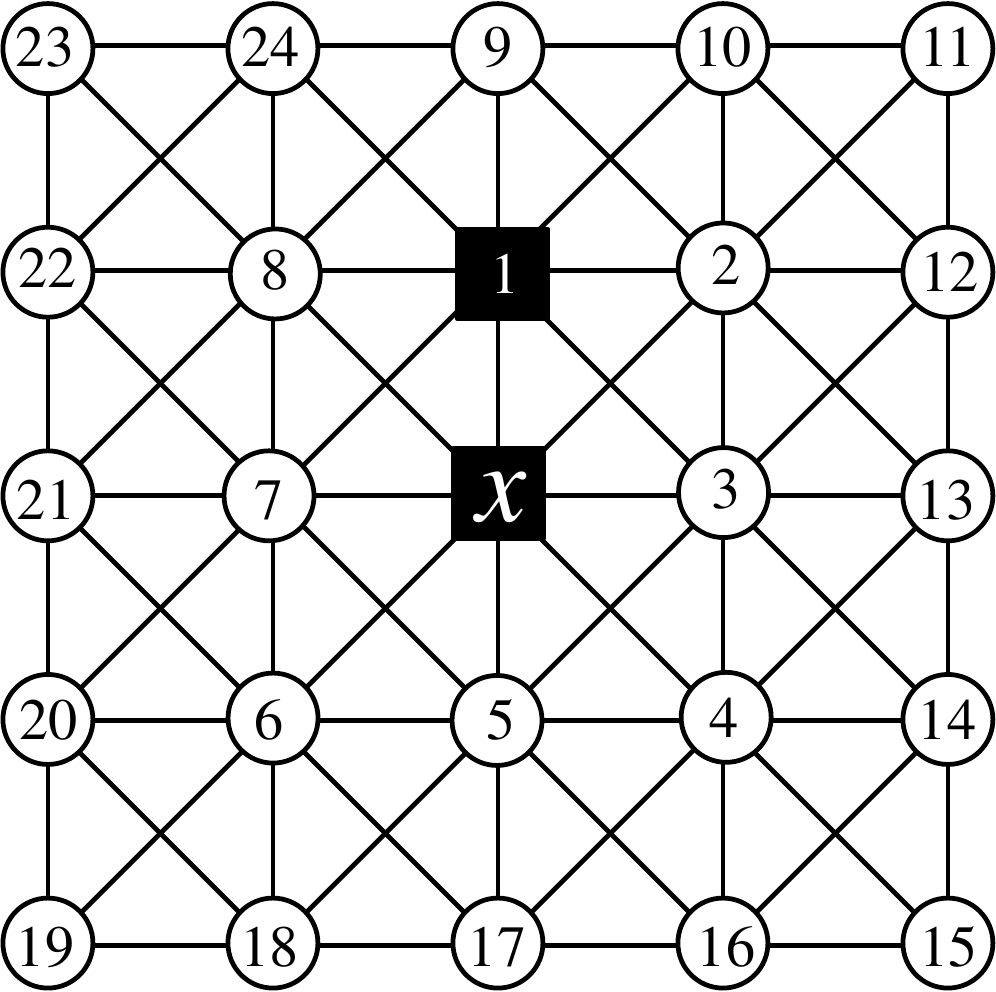}
    \caption{``Vertical''}
    \label{fig:vertical_red_ld}
\end{wrapfigure}
\newcase{Vertical}
Consider the ``vertical'' configuration given in Figure~\ref{fig:vertical_red_ld}.
We begin by examining the vertex pair $(v_2,v_3)$.
If $|\{v_2,v_3\} \cap S| = 2$ then $\{x,v_1,v_2,v_3\} \subseteq D_{4+}$, so $sh(x) \le \frac{4}{4} + \frac{5}{2} < \frac{11}{3}$ and we are done.
Next, suppose $|\{v_2,v_3\} \cap S| = 1$, then $\{x,v_1,v_2,v_3\} \subseteq D_{3+}$.
If $\{v_7,v_8\} \cap S \neq \varnothing$ then $sh(x) \le \frac{2}{4} + \frac{2}{3} + \frac{5}{2} = \frac{11}{3}$ and we are done, so we assume $\{v_7,v_8\} \cap S = \varnothing$.
We observe that if $\{v_7,v_8\} \subseteq D_2$ then they will not be distinguished, so $dom(v_7) \ge 3$ or $dom(v_8) \ge 3$ and we are done.
Lastly, suppose $|\{v_2,v_3\} \cap S| = 0$; by symmetry we can assume that $\{v_7,v_8\} \cap S = \varnothing$ as well.
We require $sh[v_2v_3] \le \max\{\sigma_{33},\sigma_{24}\}$ to be distinguished; by symmetry, this is also true for $sh[v_7v_8]$.
If $dom(v_2) = 2$, then $\{x,v_1\} \subseteq D_{3+}$ and $dom(v_3) \ge 4$ to distinguish from $v_2$; thus $sh(x) \le \sigma_{24} + \sigma_{24} + \sigma_{33} + \sigma_{222} = \frac{11}{3}$ and we are done.
Similarly, we are also done if $dom(v_3) = 2$.
Therefore, we can assume $\{v_2,v_3,v_7,v_8\} \subseteq D_{3+}$.
If $dom(x) \ge 3$ we are done, so $\{v_4,v_5,v_6\} \cap S = \varnothing$.
If $\{v_4,v_5,v_6\} \cap D_{3+} \neq \varnothing$ then we will be done, so we assume $\{v_4,v_5,v_6\} \subseteq D_2$.
Vertex $v_5$ must be 2-dominated by exactly one of $v_{16}$, $v_{17}$, or $v_{18}$.
If $v_{16} \in S$ or $v_{17} \in S$ then $v_4$ and $v_5$ will not be distinguished, a contradiction; otherwise $v_{18} \in S$, but then $v_5$ and $v_6$ are not distinguished, a contradiction.

\begin{wrapfigure}{r}{0.2\textwidth}
    \centering
    \vspace{-1em}
    \includegraphics[width=0.2\textwidth]{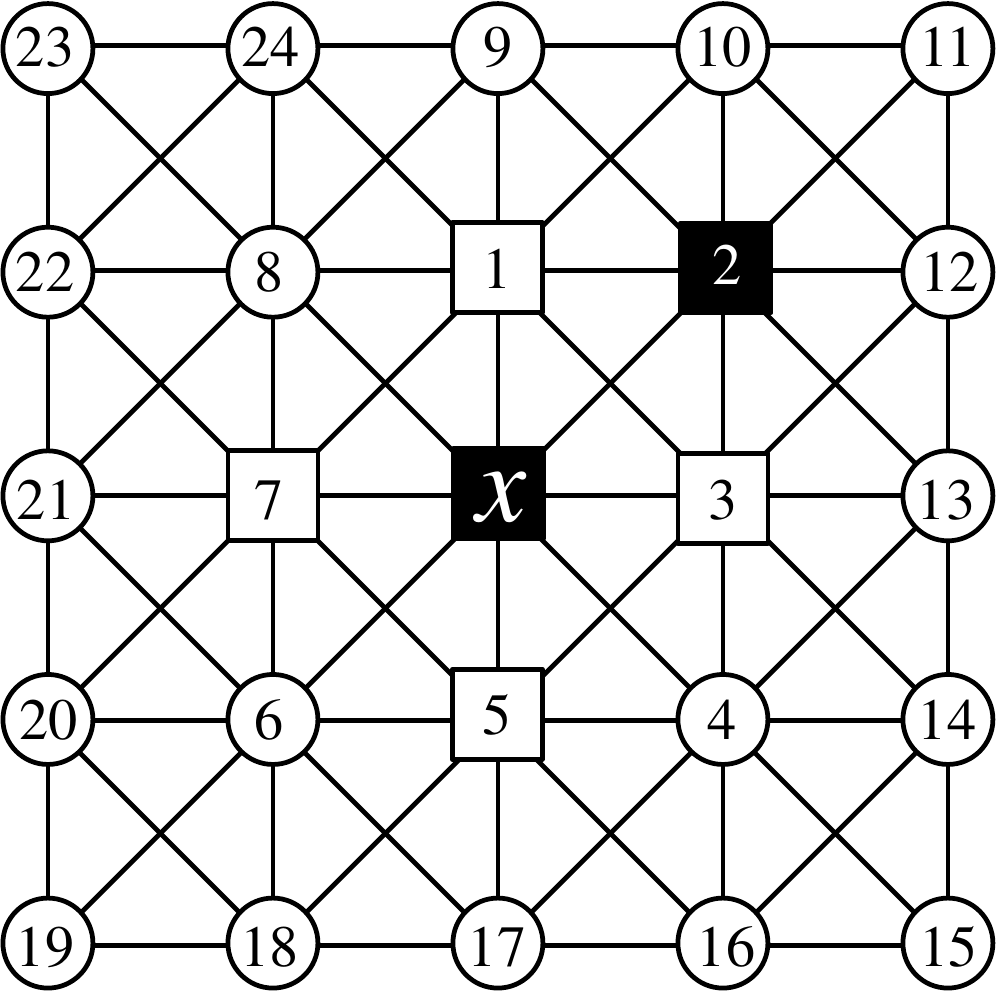}
    \caption{``Diagonal''}
    \label{fig:diagonal_red_ld}
\end{wrapfigure}

\newcase{Diagonal}
Consider the other case for 2-dominating $x$, the ``diagonal'' configuration show in Figure~\ref{fig:diagonal_red_ld}.
We can assume that $\{v_1,v_3,v_5,v_7\} \cap S = \varnothing$, as otherwise we fall into the previous case.

We observe that $sh[v_1v_3] \le \max\{\sigma_{33},\sigma_{24}\}$ to be distinguished.
Suppose $dom(v_1) = 2$, then $dom(v_3) \ge 4$ to be distinguished.
Additionally, $\{x,v_2\} \subseteq D_{3+}$ to be distinguished from $v_1$.
If $v_6 \in S$ then $sh[v_5v_7] \le \sigma_{24}$ and we are done; so we assume $v_6 \notin S$.
Suppose $v_8 \in S$, then $dom(v_7) \ge 4$ to distinguish from $v_1$ and we are done; so we assume $v_8 \notin S$.
Thus, $v_4 \in S$ to 3-dominate $x$.
If $\{v_4,v_5\} \subseteq D_2$ then they are not distinguished, so $sh[v_4v_5] \le \sigma_{23}$ and we are done.
Therefore, we can assume that $\{v_1,v_3\} \subseteq D_{3+}$.
If $v_6 \in S$ then we can assume $\{v_5,v_7\} \subseteq D_{3+}$, and $dom(x) \ge 3$, so we are done; thus, we assume $v_6 \notin S$.
Next, we consider the vertex pair $(v_4,v_8)$.
Suppose $\{v_4,v_8\} \subseteq S$.
If $\{v_4,v_5\} \subseteq D_2$ they will not be distinguished, so $sh[v_4v_5] \le \sigma_{23}$, and by symmetry $sh[v_7v_8] \le \sigma_{23}$, and we are done.
Next, we assume $|\{v_4,v_8\} \cap S| = 1$, and without loss of generality let $v_4 \in S$.
We see that $sh[v_4v_5] \le \sigma_{23}$ to be distinguished.
If $\{v_6,v_7,v_8\} \cap D_{3+} \neq \varnothing$ we will be done, so we assume $\{v_6,v_7,v_8\} \subseteq D_2$.
Any way to 2-dominate $v_7$ will result in either $(v_7,v_8)$ or $(v_6,v_7)$ not being distinguished, a contradiction.

Finally, we assume that $\{v_4,v_8\} \cap S = \varnothing$.
Suppose $v_{17} \in S$ then $dom(v_4) \ge 3$ and $dom(v_6) \ge 3$ to distinguish from $v_5$.
If $dom(v_5) \ge 3$ we are done, so we assume $dom(v_5) = 2$, so $\{v_{16},v_{17}\} \cap S = \varnothing$.
Thus, $dom(v_4) \ge 4$ and $dom(v_6) \ge 4$ to distinguish from $v_5$, and we are done; so we assume $v_{17} \notin S$, and by symmetry $v_{21} \notin S$.
If $dom(v_5) \ge 3$, then $\{v_{16},v_{18}\} \subseteq S$; we require $dom(v_4) \ge 3$ and $dom(v_6) \ge 3$ to distinguish from $v_5$, and we are done, so we can assume $dom(v_5) = 2$, and by symmetry $dom(v_7) = 2$.
Suppose $v_{16} \in S$, then $v_{18} \notin S$.
We require $dom(v_4) \ge 4$ to distinguish from $v_5$.
If $v_{22} \in S$ as well, then $v_{20} \notin S$ and we find that $dom(v_8) \ge 4$ to distinguish from $v_7$, and we are done; thus, we assume $v_{22} \notin S$, so $v_{20} \in S$ to 2-dominate $v_7$.
We find that $v_6$ and $v_7$ are not distinguished, a contradiction; therefore, we assume $v_{16} \notin S$, and by symmetry $v_{22} \notin S$.
Thus, $\{v_{18}, v_{20}\} \subseteq S$, and $dom(v_6) = 4$ with $v_{19} \in S$ to distinguish from $v_5$.
Suppose $\{v_{12},v_{13}\} \cap S \neq \varnothing$; then $dom(v_2) \ge 3$.
If $dom(v_2) \ge 4$ we are done, so we assume $dom(v_2) = 3$; thus, $\{v_9,v_{10},v_{11}\} \cap S = \varnothing$.
We require $v_{24} \in S$ to 3-dominate $v_1$ and $dom(v_8) \ge 3$ to distinguish from $v_1$, so we are done.
Therefore, we can assume $\{v_{12},v_{13}\} \cap S = \varnothing$.
Then $v_{14} \in S$ to 3-dominate $v_3$ and $dom(v_4) \ge 3$ to distinguish from $v_3$, and by symmetry we find $dom(v_8) \ge 3$ as well, completing the proof of Theorem~\ref{theo:red-sh}.
\cendproof

\begin{corollary}\label{cor:red}
$\frac{3}{11} \le \textrm{RED:LD}(\textrm{K})$
\end{corollary}

Proof of Corollary~\ref{cor:red} follows directly from Theorem~\ref{theo:red-sh}; thus, $\frac{3}{11} \le \textrm{RED:LD}(\textrm{K}) \le \frac{5}{16}$.

\begin{corollary}\label{cor:det}
$\frac{3}{11} \le \textrm{DET:LD}(\textrm{K})$
\end{corollary}

Proof of Corollary~\ref{cor:det} follows from the observation that by Theorems \ref{theo:red-ld-char} and \ref{theo:det-ld-char}, every DET:LD set is a RED:LD set \cite{jean23a}.
Thus, $\frac{3}{11} \le \textrm{DET:LD}(\textrm{K}) \le \frac{3}{8}$.

\begin{corollary}\label{cor:sh-err-ld-king}
$\frac{15}{38} \le \textrm{ERR:LD}(\textrm{K})$
\end{corollary}

Proof of Corollary~\ref{cor:sh-err-ld-king} follows directly from Theorem~\ref{theo:sh-err-ld-king}; thus, $\frac{15}{38} \le \textrm{ERR:LD}(\textrm{K}) \le \frac{7}{16}$.

\begin{conjecture}\label{conj:det}
$\frac{15}{44} \le \textrm{DET:LD}(\textrm{K})$
\end{conjecture}

We believe Conjecture~\ref{conj:det} can be proven using similar techniques and the discharging strategy introduced in Subsection~\ref{sec:discharge}.
It may be possible to prove higher lower bounds for $\textrm{RED:LD}(\textrm{K})$ and $\textrm{ERR:LD}(\textrm{K})$.
We are currently exploring other possible approaches, but we conjecture that any higher value must use either more advanced discharging techniques or expand further from the center vertex, $x$.

\bibliographystyle{ACM-Reference-Format}
\bibliography{lob-refs, ext-refs}

\end{document}